\documentclass[a4paper,10pt]{article}
\usepackage{graphicx}
\usepackage{subfigure}
\usepackage{setspace}
\usepackage{epstopdf}
\usepackage{booktabs} 
\usepackage{multirow} 
\usepackage{array}
\usepackage{threeparttable}
\usepackage{amssymb, amsmath, amsthm, epsfig}
\numberwithin{equation}{section}
\usepackage{latexsym, enumerate}
\usepackage{eepic}
\usepackage{epic}
\usepackage{color}

\usepackage{amsmath,textcomp, amssymb, mathrsfs, bm, amsthm,amsfonts}
\usepackage{algorithm}
\usepackage{algpseudocode}
\allowdisplaybreaks[4]
\usepackage{ifpdf}
\usepackage{chngcntr}
\counterwithin{figure}{section}
\usepackage{dsfont}
\usepackage{multirow}
\usepackage{bm}
\usepackage{caption}
\usepackage{float}
\usepackage[colorlinks,
            linkcolor=blue,
            anchorcolor=blue,
            citecolor=blue]{hyperref}

\theoremstyle{plain}
\newtheorem{lem}{Lemma}[section]

\newtheorem{prop}[lem]{Proposition}

\theoremstyle{definition}

\theoremstyle{remark}

\renewcommand{\thefigure}{\thesection.\arabic{figure}}

\begin{document}

\renewcommand{\figurename}{Figure}
\renewcommand{\thesubfigure}{(\alph{subfigure})}
\renewcommand{\thesubtable}{(\alph{subtable})}
\makeatletter
\renewcommand{\p@subfigure}{\thefigure~}

\makeatother
\title{\large\bf Solving Bayesian inverse problems via Fisher adaptive Metropolis adjusted Langevin algorithm}
\author{
Li-Li Wang\thanks
{School of Mathematics, Hunan University, Changsha 410082, China.
Email: lilywang@hnu.edu.cn}
\and
Guang-Hui Zheng\thanks
{School of Mathematics, Hunan University, Changsha 410082, China. Email: zhenggh2012@hnu.edu.cn (Corresponding author)}
}
\date{}
\maketitle

\begin{center}{\bf ABSTRACT}
\end{center}\smallskip

The preconditioned Metropolis adjusted Langevin algorithm (MALA) is a widely used method in statistical applications, where the choice of the preconditioning matrix plays a critical role. Recently, Titsias \cite{Titsias2024} demonstrated that the inverse Fisher information matrix is the optimal preconditioner by minimizing the expected squared jump distance and proposed an adaptive scheme to estimate the Fisher matrix using the sampling history. In this paper, we apply the Fisher adaptive Metropolis adjusted Langevin algorithm (MALA) to Bayesian inverse problems. Moreover, we provide a rigorous convergence rate analysis for the adaptive scheme used to estimate the Fisher matrix. To evaluate its performance, we use this algorithm to sample from posterior distributions in several Bayesian inverse problems. And compare its constructions with the standard adaptive Metropolis adjusted Langevin algorithm (which employs the empirical covariance matrix of the posterior distribution as the preconditioner) and the preconditioned Crank-Nicolson (pCN) algorithm. Our numerical results demonstrate show that the Fisher adaptive MALA is highly effective for Bayesian inversion, and significantly outperforms other sampling methods, particularly in high-dimensional settings.

\smallskip
{\bf keywords}: Bayesian inverse problems, Fisher adaptive, Metropolis adjusted Langevin algorithm, convergence rate

\section{Introduction}\label{sec1}
Inverse problems involve inferring unknown parameters or states of a system from observed data. These problems arise in various fields, including physics, engineering, medical imaging, and geophysics. Classical methods, such as regularization methods \cite{Tikhonov1963, Vogel2002}, provide point estimates of the unknowns but fail to quantify their uncertainty. In contrast, the Bayesian approach \cite{Kaipio2006, Stuart2010} has gained significant attention due to its ability to quantify uncertainty in  inverse problems.
Bayesian inversion applies Bayes' theorem to combine prior knowledge with the likelihood function derived from observed data, yielding a posterior distribution that includes all available information about the unknown parameters. To extract information from the posterior distribution, Markov Chain Monte Carlo (MCMC) algorithms are commonly employed. These methods provide practical tools for sampling from posterior distributions and estimating unknowns. However, traditional MCMC methods often suffer from slow convergence and high correlations, especially in high-dimensional parameter spaces or complex models.
Consequently, improving the convergence speed and accuracy of MCMC-based solutions for Bayesian inverse problems remains a critical and ongoing area of research.


The efficiency of MCMC algorithms is heavily influenced by the choice of the proposal distribution, which determines the mixing time of the Markov chain. Consequently, the design of the Markov transition kernel is crucial for achieving faster convergence and better sampling efficiency \cite{Hastings1970, Roberts2001}. Standard MCMC algorithms, such as the Random Walk Metropolis (RWM), can become arbitrarily slow as the discretization dimensionality increasing \cite{Roberts2001, Mattingly2012}.
To address this issues, a family of dimension-independent MCMC algorithms, known as the preconditioned Crank-Nicolson (pCN) algorithms, have been developed in \cite{Cotter2013}. These algorithms are derived from a Crank-Nicolson discretization of Langevin dynamics that preserves the prior measure, making them particularly robust for dimension increasing in Bayesian inverse problems. However, despite their robustness, pCN algorithms can suffer from high sample correlations and slower mixing, limiting their ability to efficiently explore complex posterior distributions.
Gradient-based methods have been proposed as an improvement to the pCN algorithm. There are a class of Markov Chain Monte Carlo (MCMC) techniques that utilize Langevin dynamics to generate samples. These methods accelerate sampling by using the local geometry of the target distribution, making them particularly effective for high-dimensional problems. Examples include Langevin Monte Carlo (LMC), the Metropolis-adjusted Langevin algorithm (MALA), stochastic gradient Langevin dynamics (SGLD) and their variants, such as \cite{Durmus2017,Durmus2019,Pasarica2016,Roberts1996,Welling2011,Zhang2023}.
Nevertheless, efficiently tuning these algorithms for posterior distributions remains a key area of research.

To further enhance MCMC performance, adaptive MCMC algorithms have been introduced. These methods automatically tune or learn the scales of the proposal distribution during sampling utilizing the chain's history, enabling the algorithm to adapt to unknown posterior structures. Adaptive algorithms improve efficiency and convergence, such as optimal scaling results in \cite{Andrieu2001, Andrieu2008, Roberts2009}. 
Recent advancements \cite{Haario2001, Rosenthal2011} have further motivated the development of adaptive methods for exploring complex posterior distributions more effectively.

Another approach to improving sampling efficiency involves the use of preconditioning matrix Metropolis adjusted Langevin algorithm (MALA). MALA is a widely used gradient-based MCMC method, derived from a first-order discretization of the continuous-time Langevin diffusion.
However, it is unclear how this matrix should be defined in any systematic and principled manner. In \cite{Titsias2024}, Titsias introduced a powerful adaptive MCMC algorithm specifically designed for high-dimensional problems. This method learns a preconditioning matrix for MALA by optimizing the expected squared jump distance \cite{Pasarica2010}. Titsias demonstrated that the optimal preconditioner is the inverse Fisher information matrix,
which minimizes the expected squared jump distance and provides an optimal scaling for the proposal distribution. This findings challenges the conventional intuition in adaptive MCMC that the covariance of the target distribution is the best preconditioner.
%
%
Building on this result, Titsias proposed an adaptive MALA algorithm that empirically estimates the inverse Fisher information matrix using the gradient history from the Markov chain. The algorithm employs a stochastic approximation scheme to recursively update the empirical Fisher matrix and adaptively tunes the step size to achieve an optimal acceptance rate of approximately $0.574$.
Numerical experiments, sampling form a given target function, demonstrated that this approach outperforms other samplers, including standard MALA (without preconditioning), covariance-based adaptive MCMC methods \cite{Haario2001}, Riemannian manifold MALA \cite{Girolami2011}, and Hamiltonian Monte Carlo (HMC).

In this work, we analyze the convergence rate of the adaptive scheme for the Fisher information matrix and apply this algorithm to Bayesian inverse problems. These aspects were not addressed in Titsias's work. Using the stochastic approximation framework, the recursive estimation of the Fisher information matrix can be extended to a more general scheme. Then provide an convergence rate analysis. Applying ideas from stochastic approximation methods presented in \cite{Andrieu2005, Kallianpur1954, Nemirovski2008}, we prove that with a learning rate of
$1/n$, where $n$ is the iterations, the convergence rate of this general scheme is $\mathcal{O}(1/n)$. Then we apply this algorithm to various Bayesian inverse problems, including both linear and nonlinear cases, and demonstrate its effectiveness through extensive numerical results.
The primary contributions of this paper are summarized as follows:
\begin{itemize}
\item
We extend the Fisher adaptive MALA algorithm to the Bayesian inversion framework.
This method uses the Fisher information matrix to construct a more efficient proposal distribution, improving sampling performance in challenging posterior distributions.
\item
We provide a rigorous theoretical analysis of the convergence properties of the adaptive Fisher matrix estimation scheme. Using a stochastic approximation framework, we prove that the proposed algorithm achieves a convergence rate of $\mathcal{O}(1/n)$, where $n$ is the iterations. 
\item
We apply this algorithm to various inverse problems, including both linear and nonlinear settings. Through extensive numerical experiments, we compare its performance with standard adaptive MALA (which uses adaptive covariance) and the pCN algorithm. The results demonstrate that the Fisher adaptive MALA outperforms the other methods in terms of convergence speed and accuracy, highlighting its practical advantages in solving Bayesian inverse problems.
\end{itemize}

The paper is organized as follows:
In section 2, we introduce Bayesian inversion framework and the Fisher information adaptive MALA. We outline the finite-dimensional Bayesian approach, describe the Fisher adaptive MALA process, and develop the convergence analysis.
Section 3 presents numerical experiments for Bayesian inverse problems that demonstrate the effectiveness of this algorithm and compare its performance with other methods.
Section 4 concludes the paper with a summary of the main findings and a discussion of potential directions for future research.

\section{Fisher Adaptive MALA in Bayesian inversion}\label{sec2}
In this section, we review the Bayesian approach to inverse problems. Then, we apply the Fisher adaptive MALA to sample from the Bayesian posterior distributions arising in inverse problems.

\subsection{Bayesian inversion framework}\label{baye}
Throughout this work, we consider the inverse problem of finding an unknown parameter $x\in \mathbb{R}^d$ from data $y\in \mathbb{R}^n$,
where the relationship between $x$ and $y$ is described by the following model
\begin{eqnarray}
\label{bayesian model}
y=\mathcal{F}(x)+\eta,
\end{eqnarray}
where $\mathcal{F}: \mathbb{R}^d \to \mathbb{R}^{n}$ is referred as forward operator, and $\eta$ represents the measurement noise, which is assumed to follow an $n$-dimensional Gaussian distribution.
In the Bayesian framework, the vectors $x$, $\eta$ and $y$ are treated as random variables. Assuming that $x$ and $\eta$ are independent with $\pi_0$ as the prior distribution and $\eta\thicksim N(0,\Sigma)$ as the Gaussian noise.
The joint distribution of $(x, y)$ is expressed as
$$\pi(x, y)=\pi(y\mid x)\pi_0(x).$$
Furthermore, applying Bayes' rule, the posterior distribution of the unknown parameter $x$ from the observed data $y$ is given by
$$\pi(x \mid y)=\frac{\pi(y| x)\pi_0(x)}{Z(y)}\propto \exp(-\Phi(x))\pi_0(x),$$
where, $Z(y)$ is the normalization constant and
\begin{equation}\label{potentialfuction}
\Phi(x) := \frac{1}{2}\big\|\mathcal{F}(x)-y\big\|_{\Sigma}^2=\frac{1}{2}(\mathcal{F}(x)-y)^{\top}\Sigma^{-1}(\mathcal{F}(x)-y),
\end{equation}
is the data fidelity term. For the prior distribution $\pi_0(x)$, we primarily assume a Gaussian prior with zero mean and covariance $C$,  i.e. $\pi_0(x)= N(0,C)$ and
\[\pi_0(x)\propto \exp\{-\frac12\|x\|^2_{C}\} =\exp\{-\frac12 x^{\top}C^{-1}x\}.\]
Thus the posterior distribution can be expressed as:
\begin{equation}\label{posterior}
\pi(x \mid y)\propto  \exp\{-\Phi(x)-\frac12\big\|x\big\|^2_C\}.
\end{equation}
Let $\pi(x)=\pi(x \mid y)$ (here omit the conditioning of data $y$), then the gradient of the log posterior is given by
\begin{equation}\label{gradposterior}
\nabla\log \pi(x)=- C^{-1}x-\nabla\Phi(x) ~\text{and} ~ \nabla\Phi(x)=(\nabla \mathcal{F}(x))^{\top}\Sigma^{-1}(\mathcal{F}(x)-y),
\end{equation}
where $\nabla \mathcal{F}(x)$ is the $Fr\acute{e}chet$ derivative of $x$.

A special case appears if $\mathcal{F}(x)$ is a linear mapping, such as $\mathcal{F}(x)=Fx$. Then it is known from \cite{Moselhy2012} that posterior distribution is Gaussian, i.e. $\pi(x\mid y)=N(\mu_{post},C_{post})$ and the mean $\mu_{post}$ and covariance matrix $C_{post}$ are given by
\begin{equation}\label{cpost}
C_{post}=(C^{-1}+F^{\top}\Sigma^{-1}F)^{-1},~\mu_{post}=C_{post}F^{\top}\Sigma^{-1}y,
\end{equation}
and
\begin{equation}\label{linear}
\nabla\Phi(x)=F^{\top}\Sigma^{-1}(Fx-y).
\end{equation}

\subsection{Fisher adaptive MALA}
In this subsection, we briefly introduce the Fisher adaptive Metropolis-Adjusted Langevin Algorithm (MALA). For further details, readers can refer \cite{Titsias2024}.

The algorithm is based on the overdamped Langevin diffusion process in continuous time
\begin{equation}
d x_t = \frac{1}{2} M \nabla \log \pi(x_t) d t + \sqrt{M} d B_t,
\label{eq:langevin}
\end{equation}
where $B_t$ is a $d$-dimensional Brownian motion with stationary distribution $\pi(x)$.
The preconditioning matrix $M$ is a symmetric positive definite covariance matrix, with $\sqrt{M}$ satisfying $\sqrt{M}^\top \sqrt{M}= M$.
The preconditioned MALA is derived by discretizing the Langevin diffusion equation \eqref{eq:langevin} using the first-order Euler-Maruyama method and combining with Metropolis-Hastings adjustment.
Given the current state $x_n$ (where $n=1,2,\ldots$), the proposal state $y_n$ is generated from the proposal distribution
\begin{equation}
q(y_n \mid x_n)  = \mathcal{N}(y_n \mid x_n + \frac{\sigma^2}{2}
M \nabla \log \pi(x_n), \sigma^2
M),
\label{eq:MALAproposal}
\end{equation}
where $\sigma^2 > 0 $  is the step size.
Following  \cite[Proposition 1]{Titsias2024}, the Metropolis acceptance probability $$\alpha(x_n, y_n) = \min\left(1, r_n \right)$$ is determined by the Metropolis ratio:
$$
r_n = \frac{\pi(y_n)}{\pi(x_n)}\frac{q(x_n\mid y_n)}{q(y_n\mid x_n)},
$$
where the transition probability ratio $\frac{q(x_n\mid y_n)}{q(y_n\mid x_n)}$ is computed as
\begin{equation}\label{hxy}
\frac{q(x_n\mid y_n)}{q(y_n\mid x_n)}
 = \exp\{
h(x_n,y_n) - h(y_n,x_n)
\}, \quad
 h(u, v)
= \frac{1}{2}
\left(u -  v -  \frac{\sigma^2}{4} M \nabla \log \pi(v)
\right)^\top \nabla \log \pi(v).
\end{equation}

To derive the optimal preconditioner, consider the rejection-free or unadjusted Langevin sampler.
Discretizing the time continuous Langevin diffusion equation in \eqref{eq:langevin} with a small finite $\delta :=\sigma^2 > 0$
so that
\begin{equation}
x_{t + \delta} - x_t = \frac{\delta}{2} M \nabla \log \pi(x_t)  + \sqrt{M} (B_{t+ \delta} - B_t),  \ \ \text{where} \ \  B_{t+ \delta} - B_t \sim \mathcal{N}(0,\delta I).
\label{eq:discretizesLangevin}
\end{equation}
The expected squared jumped distance $J(\delta, M) = E[||x_{t+\delta}  - x_t ||^2]$
is calculated  as \cite[Proposition 2]{Titsias2024},
\[E[(x_{t + \delta} - x_t) (x_{t + \delta} - x_t )^\top]
= \frac{\delta^2}{4}
M E_{\pi(x_t)}
\left[ \nabla \log \pi(x_t)
\nabla \log \pi(x_t)^\top \right] M
+ \delta M,\]
\[\text{tr}\left( E[(x_{t + \delta} - x_t) (x_{t + \delta} - x_t )^\top] \right)
= E[ \text{tr}\left((x_{t + \delta} - x_t) (x_{t + \delta} - x_t )^\top \right)] =
  E[||x_{t + \delta} - x_t||^2].
\]

According to \cite[Proposition 3]{Titsias2024}, follow the common practice to tune parameter $\sigma^2$ to achieve an average acceptance rate around $0.574$, as suggested by optimal scaling results (see \cite{Andrieu2001, Andrieu2008, Roberts2009}). Under the trace constraint $tr(M)=c$, minimizing the expected squared jump distance $J(\delta, M)$ yields the optimal preconditioner:
\begin{equation}\label{eq:optimalA_and_fisher}
M^* \propto \mathcal{I}^{-1}, \ \  \mathcal{I} = E_{\pi(x)}
\left[ \nabla \log \pi(x)
\nabla \log \pi(x)^\top \right].
\end{equation}
where $\mathcal{I}$ is the Fisher information matrix. This result implies that the preconditioner $M$ should be proportional to the inverse Fisher information matrix $\mathcal{I}^{-1}$.
Specifically, the normalized proposal distribution is defined as
\begin{equation}
q(y_n\mid x_n) = \mathcal{N}\left(y_n|x_n + \frac{\sigma^2}{ \frac{2}{d} \text{tr}(M)}
M \nabla \log \pi(x_n), \frac{\sigma^2}{\frac{1}{d} \text{tr}(M)}
M\right),
\label{eq:normalizedPMALA}
\end{equation}
where $\sigma^2$ is normalized by the average eigenvalue of $M$, i.e.\ $\frac{1}{d}\text{tr}(M)$.

Following \cite{Titsias2024}, the preconditioner $M$ is adaptively estimated through the gradient history. Let the score function at the $n$-th MCMC iteration be $s_n := \nabla_{x_n} \log \pi(x_n)$. An empirical estimate of the Fisher information matrix is given by \begin{equation}
\hat{\mathcal{I}}_n = \frac{1}{n}
\sum_{i=1}^n s_i s_i^\top + \frac{\lambda}{n} I_d,
\label{eq:empiricalFisher}
\end{equation}
where $\lambda > 0$ is a fixed damping parameter. The preconditioner $M_n$ is then proportional to the inverse of $\mathcal{I}_n$
\begin{equation}
M_n \propto \left( \frac{1}{n} \sum_{i=1}^n s_i s_i^\top + \frac{\lambda}{n} I_d \right)^{-1}
=  n \left(\sum_{i=1}^n s_i s_i^\top + \lambda I_d \right)^{-1}.
\end{equation}
Using the Woodbury matrix identity, $M_n$ can be updated iteratively as
\begin{align}
&
\text{Initialization:} \
M_1 = \left(
s_1 s_1^\top + \lambda I_d
\right)^{-1}
= \frac{1}{\lambda}
\left(I_d - \frac{s_1 s_1^\top}{\lambda + s_1^\top s_1} \right),  \\
& \text{Iteration:} \  M_n  = \left(M_{n-1}^{-1} + s_{n} s_{n}^\top  \right)^{-1} = M_{n-1} -
\frac{M_{n-1} s_n s_n^\top M_{n-1}}{1 + s_n^\top M_{n-1} s_n}.
\end{align}
In practice, only the square root matrix $R_n := \sqrt{M}_n$ is required, such that $R_n R_n^\top = M_n$.  It is updated as
\begin{align}
& \text{Initialization:}\label{propSQRT1}
\ R_1
= \frac{1}{\sqrt{\lambda}}
\left(I_d - r_1 \frac{s_1 s_1^\top}{\lambda + s_1^\top s_1} \right),  \ \ r_1 = \frac{1}{1 + \sqrt{\frac{\lambda}{\lambda + s_1^\top s_1}}} \\
&
\text{Iteration:}\label{propSQRT2}
 \ R_n
= R_{n-1} -
r_n \frac{(R_{n-1} \phi_n) \phi_n^\top}{1 + \phi_n^\top \phi_n}, \ \ \phi_n = R_{n-1}^\top s_n, \ r_n = \frac{1}{1 + \sqrt{\frac{1}{1 + \phi_n^\top \phi_n}}}.
\end{align}

Using the stochastic approximation framework, the Fisher matrix $\hat{\mathcal{I}}$ can be updated as
\begin{equation}\label{eq:empiricalFisherGamma}
\hat{\mathcal{I}}_n =
\hat{\mathcal{I}}_{n-1} + \gamma_n (s_n s_n^\top - \hat{\mathcal{I}}_{n-1})
= (1 - \gamma_n) \hat{\mathcal{I}}_{n-1}
+ \gamma_n s_n s_n^\top,
\end{equation}
where the sequence is initialized at $\hat{\mathcal{I}}_1 = s_1 s_1^\top + \lambda I$, and $\gamma_n$ is learning rate, satisfying $\sum_{n=1}^{\infty} \gamma_n=\infty$, $\sum_{n=1}^{\infty} \gamma_n^2<\infty$.

Under the bounded fourth-moment condition $E[\|s_ns_n^\top\|^2_{F}]<\infty$, convergence of the Fisher information matrix estimate is guaranteed. 
Using the stochastic approximation results from \cite{Andrieu2005, Kallianpur1954, Nemirovski2008}, we can prove that the following proposition holds.
\begin{prop}
\label{prop:convergence}
Let $\pi(x)$  be the target distribution, and define
$s(x)=\nabla\log\pi(x)$ and $x\thicksim \pi(x)$. Denote $s_n=\nabla\log\pi(x_n)$ where $x_n\thicksim \pi(x_n)$. Assume fourth moment of the score function $E[\|s_ns_n^\top\|^2_{F}]$ is bounded. If the learning rates is set as $\gamma_n=1/n$, then the online learning update for the
empirical Fisher information matrix, as described in \eqref{eq:empiricalFisherGamma} achieves the asymptotically optimal convergence rate $\mathcal{O}(1/n)$, i.e.
\[E[\|\hat{\mathcal{I}}_{n}-\mathcal{I}\|^2_{F}]=\mathcal{O}(1/n).\]
\begin{proof}
Let $\xi=ss^\top$ and $\xi_n=s_ns_n^\top$.
Our aim is to estimate $E[\xi]$, and $\lim_{n\to \infty}\hat{\mathcal{I}}_{n}=E[\xi]=\hat{\mathcal{I}}$. Stochastic approximation to a broad class of stochastic iterative algorithm solving root finding or optimization problems. Suppose we want to find the optimization problems
\begin{equation}
\min_{z\in \mathcal{X}}\{f(z)=E[\frac12\|z-\xi\|_F^2]\}, ~\mathcal{X}=R^{d\times d}.
\end{equation}
Thus we would like to find the root of the equation
\begin{equation}
\nabla f(z)=E[z-\xi]=0.
\end{equation}
Define $g(z,\xi)=z-\xi$, which is an unbiased estimate of $\nabla f$ and referred as the subgradient.
Therefore, the stochastic approximation framework \eqref{eq:empiricalFisherGamma} (where ${\mathcal{I}}_{n}$ is iter variable of $z$)
can solve the above problem, i.e.
$$\begin{aligned}\hat{\mathcal{I}}_n
&= \hat{\mathcal{I}}_{n-1}-\gamma_n g(\hat{\mathcal{I}}_{n-1},\xi_n)\\
&= \hat{\mathcal{I}}_{n-1}-\gamma_n (\hat{\mathcal{I}}_{n-1}-\xi_n)
\end{aligned}$$
Denote
 $ a_{n}=E[\|\hat{\mathcal{I}}_{n}-\hat{\mathcal{I}}\|_F^2],$
then we can write
\begin{equation}
\label{eq:squarederror2}
\begin{aligned}
a_{n}&=E[\|\hat{\mathcal{I}}_{n-1} -\gamma_n (\hat{\mathcal{I}}_{n-1}-\xi_n)-\hat{\mathcal{I}}\|_F^2]\\
&=E[\|\hat{\mathcal{I}}_{n-1}-\hat{\mathcal{I}}\|_F^2]+\gamma_n^2 E[\|\hat{\mathcal{I}}_{n-1}-\xi_n\|_F^2
-2\gamma_nE[({\mathcal{I}}_{n-1}-\hat{\mathcal{I}})^\top(\hat{\mathcal{I}}_{n-1}-\xi_n)].
\end{aligned}
\end{equation}
Since $\hat{\mathcal{I}}_{n-1}$ depends on $\xi_1,\xi_2,\cdots,\xi_{n-1}$ but not on $\xi_n$, we obtain
$$
\begin{aligned}
E[({\mathcal{I}}_{n-1}-\hat{\mathcal{I}})^\top(\hat{\mathcal{I}}_{n-1}-\xi_n)]
&=E\{E_{\xi_n}[({\mathcal{I}}_{n-1}-\hat{\mathcal{I}})^\top(\hat{\mathcal{I}}_{n-1}-\xi_n)|\xi_1,\xi_2,\cdots,\xi_{n-1}]\}\\
&=E\{({\mathcal{I}}_{n-1}-\hat{\mathcal{I}})^\top E_{\xi_n}[(\hat{\mathcal{I}}_{n-1}-\xi_n)|\xi_1,\xi_2,\cdots,\xi_{n-1}]\}\\
&=E[({\mathcal{I}}_{n-1}-\hat{\mathcal{I}})^\top\nabla f({\mathcal{I}}_{n-1},\xi_n)].
\end{aligned}
$$
Noting that $\nabla^2 f=I$, the expectation function $f(z)$ is strongly convex.
Thus the minimizer ${\mathcal{I}}$ is unique and the following is holds
$$({\mathcal{I}}_{n-1}-\hat{\mathcal{I}})^\top\nabla f({\mathcal{I}}_{n-1},\xi_n)\geq\|\hat{\mathcal{I}}_{n-1}-\hat{\mathcal{I}}\|_F^2.$$
By the assumption $E[\|s_ns_n^\top\|^2_{F}]=E(\xi_n\xi_n^\top)$ is bounded, then
$$\begin{aligned}
E[\|\xi_n-\hat{\mathcal{I}}_{n-1}\|_F^2]&\leq 2E[\xi_n^\top\xi_n]
+2E[\hat{\mathcal{I}}_{n-1}^\top\hat{\mathcal{I}}_{n-1}]\\
&\leq M^2.
\end{aligned}$$
Therefore, it follows that by taking expectation of both sides of \eqref{eq:squarederror2}, we have
$$a_n \leq (1-2\gamma_n)a_{n-1}+\gamma_n^2M^2.$$
If the learning rate $\gamma_n=1/n$,  we have
$$a_n \leq (1-2/n)a_{n-1}+M^2/n^2.$$
It follows by induction that
$a_n\leq Q/n$, ~where $Q=\max\{M^2,\|\hat{\mathcal{I}}_{1}-\hat{\mathcal{I}}\|_F^2\}$.
When $n=1$, $a_1=\|\hat{\mathcal{I}}_{1}-\hat{\mathcal{I}}\|_F^2$, it obviously holds. Suppose $a_{n-1}$ the result holds, therefore
$$\begin{aligned}
a_n&\leq (1-\frac2n)\frac{Q}{n-1}+\frac{M^2}{n^2}\\
&\leq\frac{(n-2)Q}{n(n-1)}+\frac{M^2}{n^2}\\
&=\frac{Q}{n}-\frac{Q}{n(n-1)}+\frac{M^2}{n^2}\\
&\leq \frac{Q}{n}.
\end{aligned}$$
Thus, we have
\[E[\|\hat{\mathcal{I}}_{n}-\hat{\mathcal{I}}\|_F^2]\leq \frac{Q}{n}.\]
The proof is completed.\\
\end{proof}
\end{prop}
Combined with the Fisher adaptive MALA, we can express the algorithm of drawing from the posterior distribution for Bayesian inversion as outlined by Algorithm \ref{alg:FisherMALA}.

\begin{algorithm}[tb]
   \caption{Fisher adaptive MALA for Bayesian inversion}
   \label{alg:FisherMALA}
\begin{algorithmic}
    \State \textbf{Input:}
        \begin{itemize}
            \item Number of samples $N_s$
            \item Damping parameter $\lambda = 10$
            \item Target acceptance rate $\alpha_* = 0.574$
            \item Step size adaptation rate $\rho_n = 0.015$
            \item Log-posterior $\log\pi(x)$ (equation \eqref{posterior})
            \item Fréchet derivative $\nabla\mathcal{F}(x)$ and gradient $\nabla\log\pi(x)$ (equation \eqref{gradposterior} or \eqref{linear})
        \end{itemize}

   \State \textbf{Initialization:}
        \begin{itemize}
            \item Run standard MALA for $n_0 = 500$ iterations to:
                \begin{itemize}
                    \item Initialize chain state $x_1$
                    \item Tune $\sigma^2$ is adapted so that acceptance rate $\alpha$ towards $\alpha_*$
                \end{itemize}
            \item Set square root matrix $R_1 = I_d$
            \item Compute normalized step size $\sigma_R^2 = \sigma^2 / \frac{1}{d}\text{tr}(R_1R_1^\top)$
            \item Evaluate $(\log\pi(x_1), \nabla\log\pi(x_1))$
        \end{itemize}
        \setstretch{1.35}
   \For{$n=1$ {\bf to} $N_{s}$}
   \State:  Propose $y_n = x_n + (\sigma_R^2 /  2) R (R^\top \nabla \log \pi(x_n))  + \sigma_R R \eta, \  \ \eta \sim \mathcal{N}(0, I_d)$
   \State: Compute $(\log \pi(y_n), \nabla \log \pi(y_n))$ use \eqref{posterior} and \eqref{gradposterior}.
   \State: Compute acceptance probability
   $\alpha(x_n, y_n) = \text{min}\left(1, e^{
\log \pi(y_n) + h(x_n,y_n) - \log \pi(x_n) - h(y_n,x_n)}
 \right)$.
   \State: {Compute adaptation signal $s_n^{\delta}  = \sqrt{\alpha(x_n, y_n)} (\nabla \log \pi(y_n) -  \nabla \log \pi(x_n))$ (this adapts by the Rao-Blackwellized score function increments $s_n^{\delta}  = \sqrt{\alpha(x_n, y_n)} (s(y_n) -  s(x_n))$) to reduce variance.}
   \State: {Use $s_n^\delta$ to adapt $R$: if $n=1$ use \eqref{propSQRT1} and if $n>1$ use \eqref{propSQRT2} (i.e. learning rate $\gamma_n=\frac{1}{n}$ in equation \eqref{eq:empiricalFisherGamma})}.
   \State: {Adapt step size $\sigma^2 \leftarrow \sigma^2 \left[1 + \rho_n (\alpha(x_n, y_n) -  \alpha_*)\right]$}.
    \State: {Normalize step size $\sigma_R^2 = \sigma^2 / \frac{1}{d} \text{tr}(R R^\top)$}.
   \State: Accept/reject $y_n$ with probability $\alpha(x_n, y_n)$ to obtain  $(x_{n+1}, \log \pi(x_{n+1}), \nabla \log \pi(x_{n+1}))$.
   \EndFor
\end{algorithmic}
\end{algorithm}

\section{Numerical Examples}\label{sec5}
In this section, we present several numerical experiments,
employing the Fisher Adaptive Metropolis adjusted Langevin algorithm (hereinafter abbreviated as FisherMala) for solving inverse problems.  For comparison, we also use two alternative samplers:
 {\bf (i)} A standard adaptive Metropolis adjusted Langevin algorithm (hereinafter abbreviated as AdaMala) where the proposal follows exactly the from in
 \eqref{eq:normalizedPMALA} but but with a preconditioning matrix that learns through standard adaptive Markov Chain Monte Carlo (MCMC).
This approach employs the recursion derived by \cite{Haario2001}, updating the mean vector $\mu_n$ and covariance matrix $\mathcal{C}n$ as follows:
\begin{equation}
\mu_n = \frac{n-1}{n}
\mu_{n-1} + \frac{1}{n} x_n,  \
\mathcal{C}_n
= \frac{n-2}{n-1}
\mathcal{C}_{n-1} + \frac{1}{n}
(x_n - \mu_{n-1}
)(x_n - \mu_{n-1})^\top,
\label{eq:adaMALA}
\end{equation}
The initial conditions are set as $\mu_1= x_1$
and
$\mathcal{C}_2 = \frac{1}{2} (x_2 - \mu_1)(x_2 - \mu_1)^\top +  \lambda I$,  and $\lambda>0$ is the damping parameter.
{\bf (ii)} A standard preconditioned Crank-Nicolson (pCN) sampler with proposal
$N(y_n\mid\sqrt{1-\beta^2} x_n, \beta^2 C)$ therein \cite{Beskos2008,Cotter2013}.
Besides, for all sampler, the initial points are drawing from prior distribution, and the initialization phase described in Appendix \ref{initial}. The Gaussian noise $\eta$ are distributed as $N(0,\epsilon^2 I)$, where the $\epsilon$ referred as noisy level.

For the linear inverse problem as Example \ref{example1}, the gradient of logarithmic posterior can be computed directly. However, when dealing with the nonlinear inverse inverse problems presented in Example \ref{example2} and Example \ref{example3}, we employ a finite difference method to approximate the $Fr\acute{e}chet$ derivatives of forward maps, as proposed in \cite{Li2008}.

\subsection{Inverse source problem}\label{example1}
In this example, we consider the following initial-boundary value problem for the non-homogeneous heat equation
\begin{equation}\label{problem1}
\begin{split}
\begin{cases}
&\dfrac{\partial u(x,t)}{\partial t}-\Delta u(x,t)=f(x), ~~~ (x,t)\in\Omega\times(0,T],\\
&u(x,t)=0, ~~~ (x,t)\in\partial\Omega\times(0,T],\\
&u(x,0)=\sin(\pi x), ~~~ x\in\Omega.\\
\end{cases}
\end{split}
\end{equation}
where the objective is to determine the heat source $f(x)$ from the final temperature measurements $u(x,T)$ for $x\in \Omega$.

Using the finite difference method (FDM) as  described in  \cite{Yan2010}, this inverse problem can  be reduced to solve the following linear inverse problem:
\begin{equation}\label{probleminverse}
y=Ff+\eta,
\end{equation}
where $\eta$ is the Gaussian observed noise. In this simulation, we fix the noise level $\epsilon = 0.01$, set the spatial domain as $\Omega = (0,1)$, and choose the time interval $T = 1$.
The exact heat source, we need to reconstruct, is given by
$$ f(x)= 2\pi^2\sin(\pi x).$$

To avoid the so-called "inverse crime",  we adopt a fine mesh for the forward problem and a coarser one for the inverse problem. For the inverse problem, the spatial domain is discretized with $N_x = 100$ and $600$ uniform grids, while the time domain uses $ N_t = 100$. Our aim is to recover the unknown function $f(x)$ in both low-dimensional ($d=100$) and high-dimensional ($d=600$) cases, and to compare the performance of each sampler.
For all three samplers, we draw $2\times10^5$ samples from posterior distribution, with the first $1\times10^5$ samples discarded as the burn-in phase. The prior distribution $\pi_0$ is assumed to have a spatially dependent covariance matrix, which defined as following:
\begin{equation}\label{cov}
c_{0}(x_1,x_2)=\gamma\exp\left[-\frac{1}{2}\left(\frac{x_1-x_2}{l}\right)^2\right].
\end{equation}
\begin{table}[htpb]
\belowrulesep=0pt
\aboverulesep=0pt
\centering
\caption{Numerical errors in Example 1.}
\setlength{\tabcolsep}{5pt} 
\renewcommand{\arraystretch}{2}
\begin{tabular}{c|ccc|ccc}
\toprule
{} &\multicolumn{3}{c|}{$d=100$}&\multicolumn{3}{c}{$d=600$}\\[-1em]
{} & {pCN} & {AdaMala} & {Fisher} & {pCN} & {AdaMala} & {Fisher}\\
\midrule
err$(\%)$ & 0.99 & 0.74 & 0.71 & 0.98 & 0.80 & 0.64 \\
\bottomrule
\end{tabular}
\label{ftab}
\end{table}
\begin{figure}[htpb]
\centering
\subfigure[d=100]{
\label{f:fig1}
\includegraphics[width=0.43\textwidth,trim=10 0 35 18,clip]{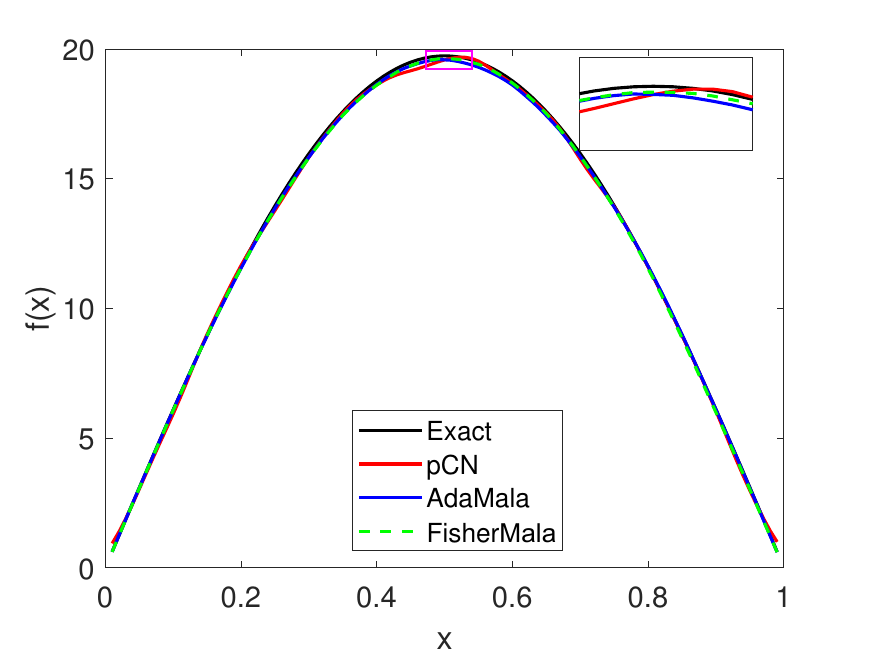}
}
\subfigure[d=600]
{
\label{f:fig2}
\includegraphics[width=0.43\textwidth,trim=10 0 35 18,clip]{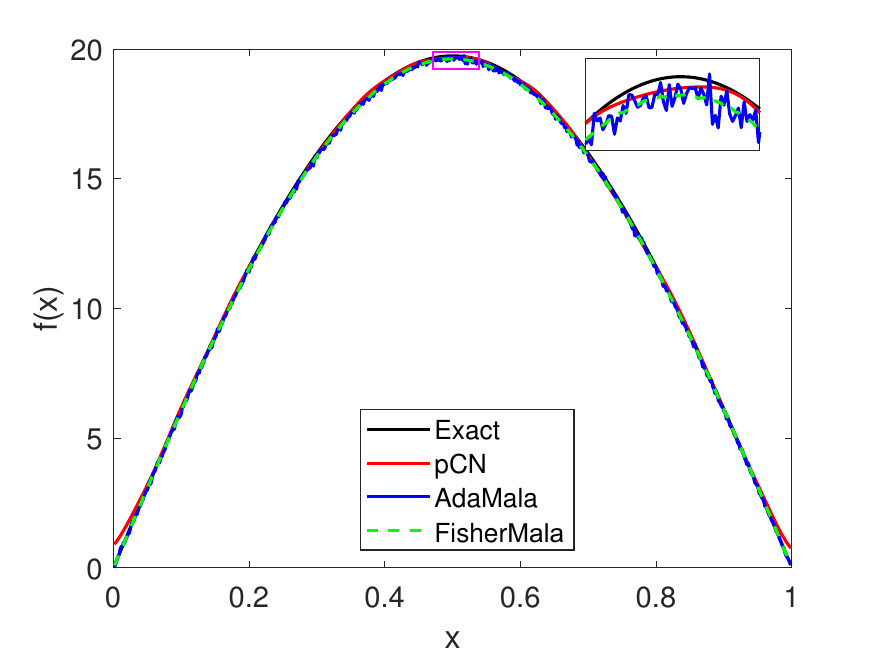}
}
\vspace{-0.3cm}
\caption{The exact solution and inversion results.
}
\label{Isf}
\end{figure}
\begin{figure}[htpb]
\centering
{
\label{f10:fig1}
\includegraphics[width=0.31\textwidth,trim=15 0 35 5,clip]{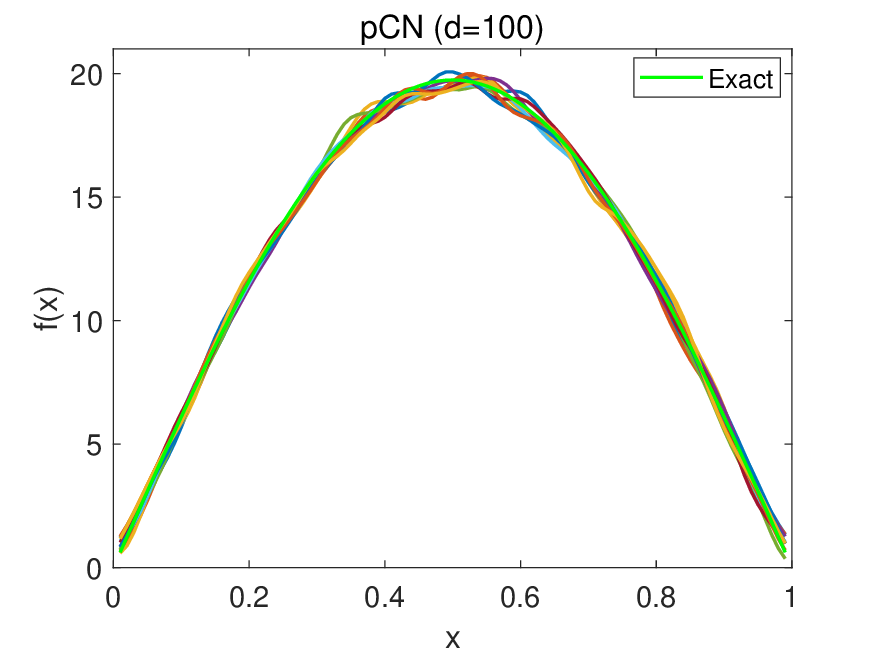}
}
{
\label{f10:fig2}
\includegraphics[width=0.31\textwidth,trim=15 0 35 5,clip]{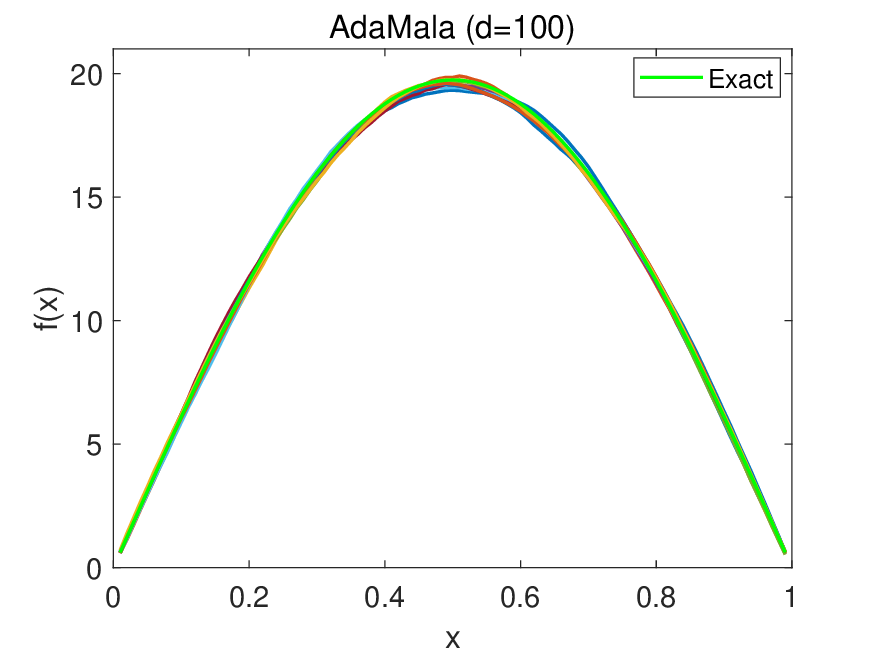}
}
{
\label{f10:fig3}
\includegraphics[width=0.31\textwidth,trim=15 0 35 5,clip]{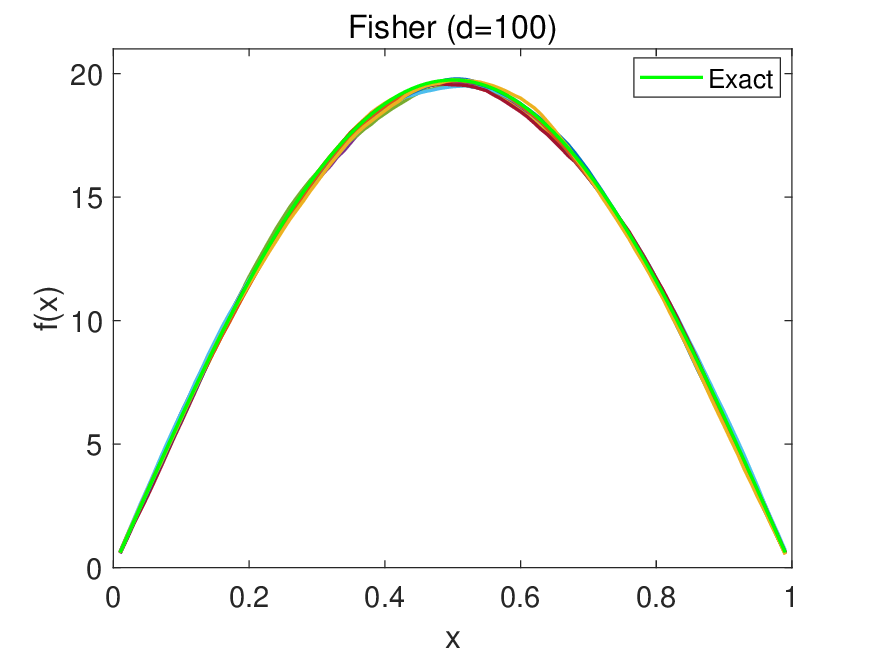}
}
\vspace{0.3cm}\\
{
\label{f10:fig4}
\includegraphics[width=0.31\textwidth,trim=15 0 40 5,clip]{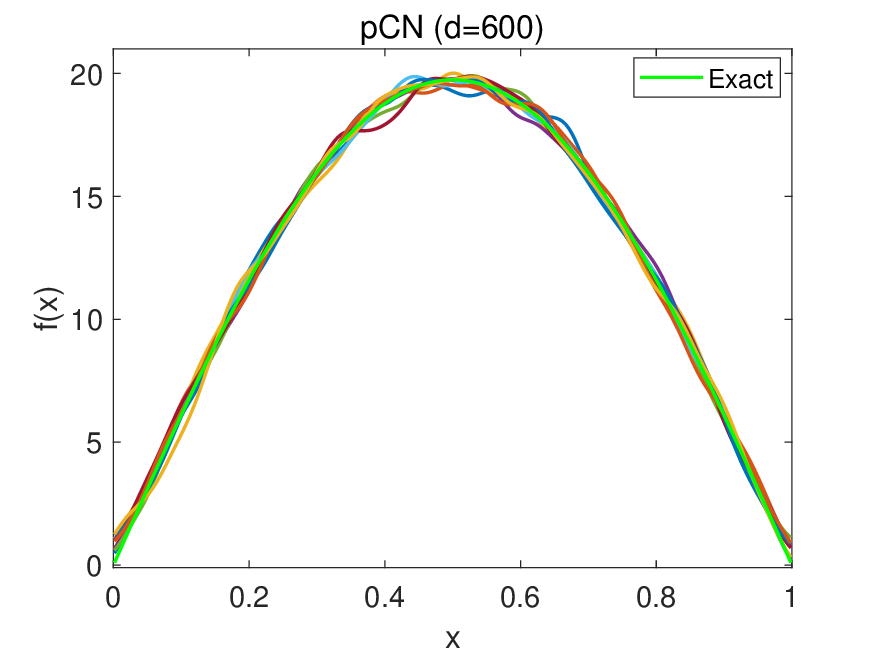}
}
{
\label{f10:fig5}
\includegraphics[width=0.31\textwidth,trim=15 0 35 5,clip]{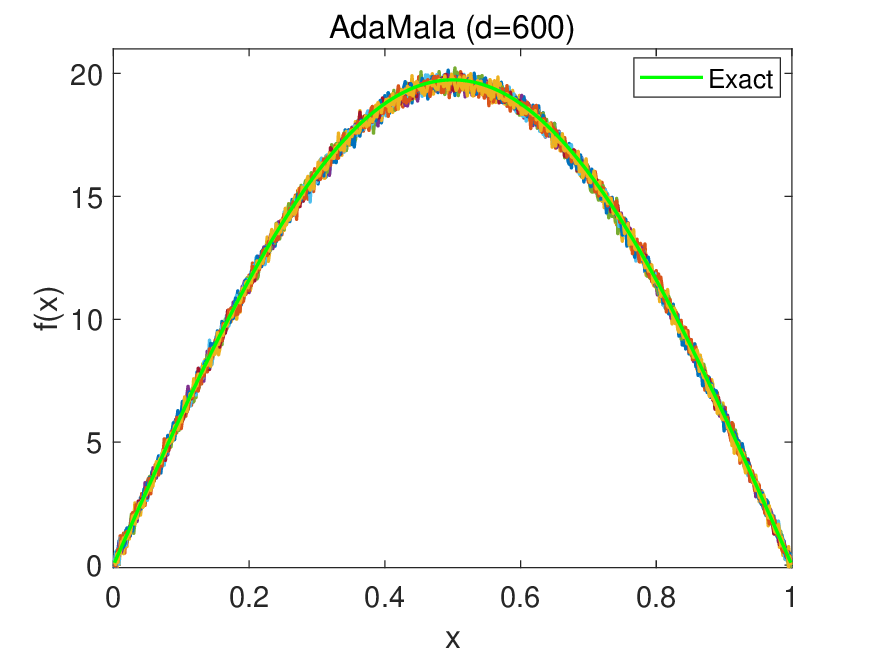}
}
{
\label{f10:fig6}
\includegraphics[width=0.31\textwidth,trim=15 0 35 5,clip]{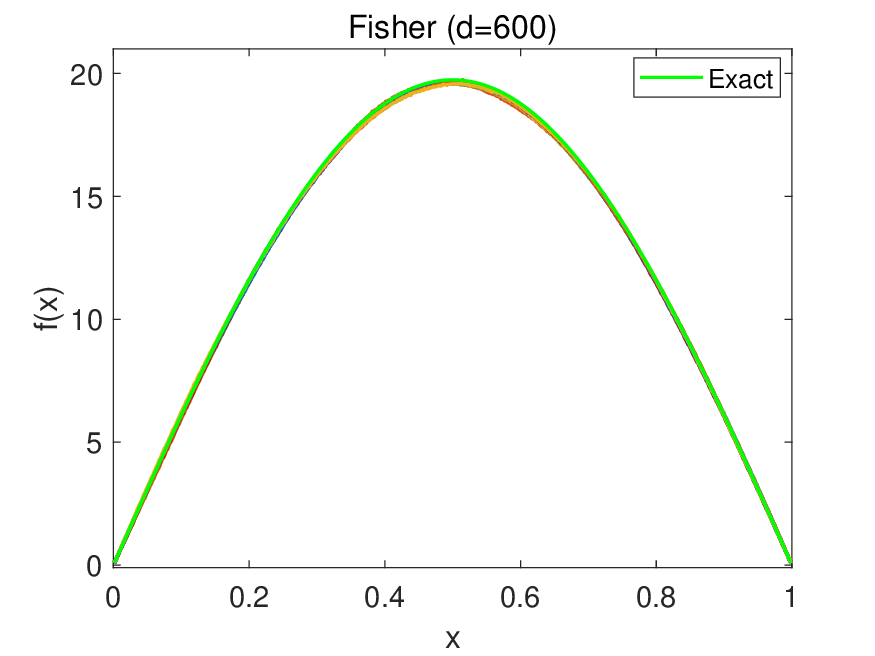}
}
\caption{The inversion results of ten-time sampling after burn-in phase for all sampler.
}
\label{Isf10}
\end{figure}

For the pCN algorithm, we set $\gamma=0.2$, $l=0.03$, and use a fixed step size $\beta=0.02$. The average acceptance rate is around $25\%$ for both $d=100$ and $d=600$, demonstrating the dimension-independent property of the pCN algorithm, as stated in \cite{Cotter2013}.
After numerous trials to tune the parameters in prior covariance matrix \eqref{cov}, we observed that the parameter $l$ must be smaller so that the prior covariance matrix to become a scalar matrix for FisherMala and AdaMala sampler. Thus, we take the prior distribution $\pi_0=N(0, 1.5 I_d)$ in FisherMala and AdaMala.
Numerical inversion results are obtained by averaging performance over ten independent runs. To assess the accuracy of the three sampling algorithms, we utilize the relative error between the exact solutions $x$ and the numerical approximations $\hat{x}$, defined by
 \[err (\%)= \dfrac{\|\hat{x}-x\|_2}{\|x\|_2}\times 100\%.\]
The reconstruction results and corresponding relative errors are presented in Figure \ref{Isf} and Table \ref{ftab}. We can see the inversion results of all three samplers are closely approximating the exact solution. However, the FisherMala performs the best, particularly in the high-dimensional $d=600$.
Figure \ref{Isf10} displays ten independent sampling realizations with varying initial values and random noise. We note that all of the numerical results can approximate the exact solutions. In the detailed information of reconstruction results, the FisherMala demonstrates robust performance, while the AdaMala and pCN samplers lead to small oscillations.
These results conclusively demonstrate that the FisherMala sampler outperforms both AdaMala and pCN, particularly in high-dimensional settings ($d=600$).
\begin{figure}[htpb]
\centering
\subfigure[d=100]{
\label{Isacf:fig1}
\includegraphics[width=0.43\textwidth,trim=10 0 25 18,clip]{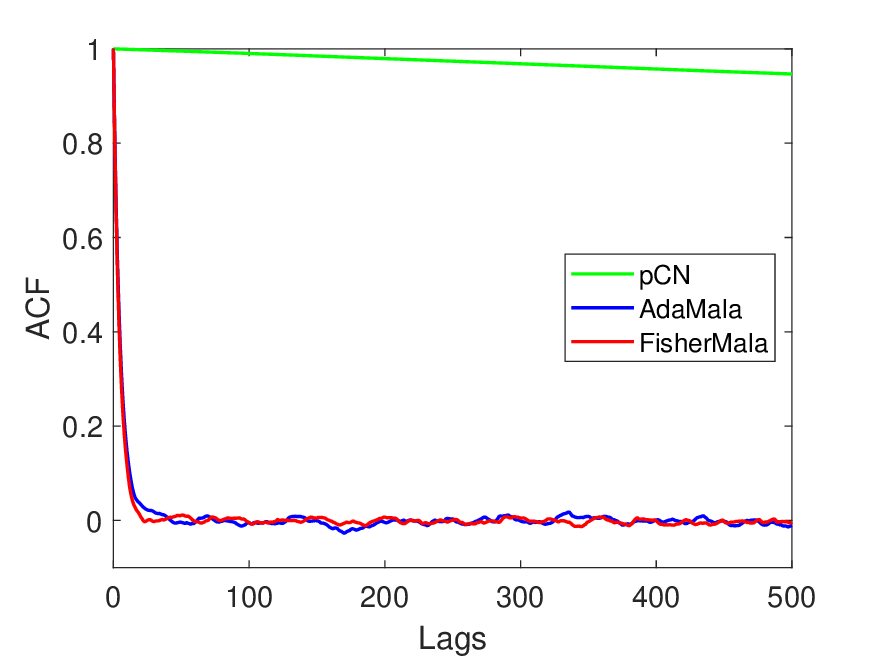}
}
\subfigure[d=600]{
\label{Isacf:fig2}
\includegraphics[width=0.43\textwidth,trim=10 0 25 18,clip]{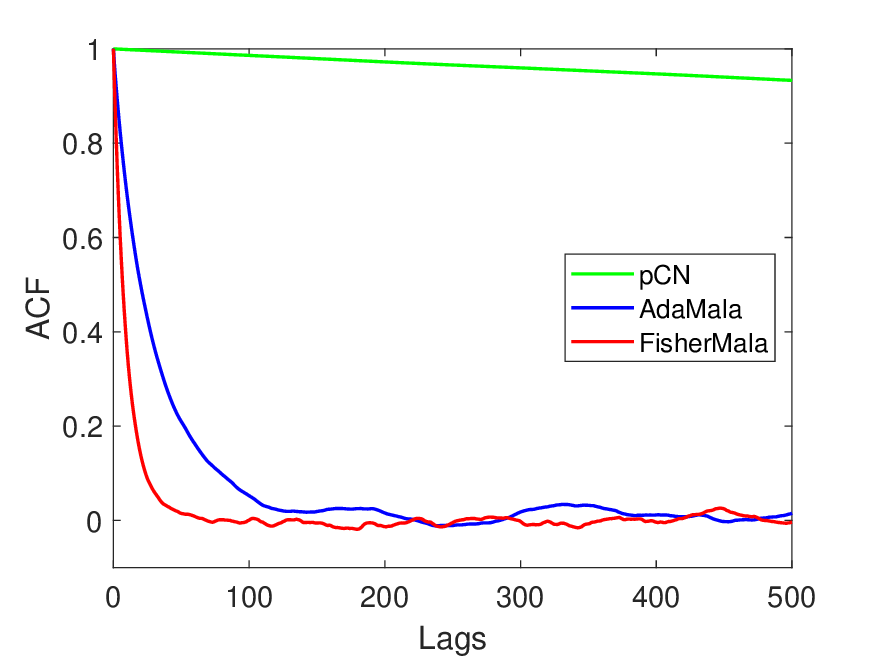}
}
\caption{The ACF for the point of state space $x_p=0.5$ for unknown function $f(x)$ with $lags=500$.}
\label{Isacf}
\end{figure}
\begin{figure}[htpb]
\centering
\subfigure[d=100]{
\label{fig2:subfig3}
\includegraphics[width=0.43\textwidth,trim=4 0 35 18,clip]{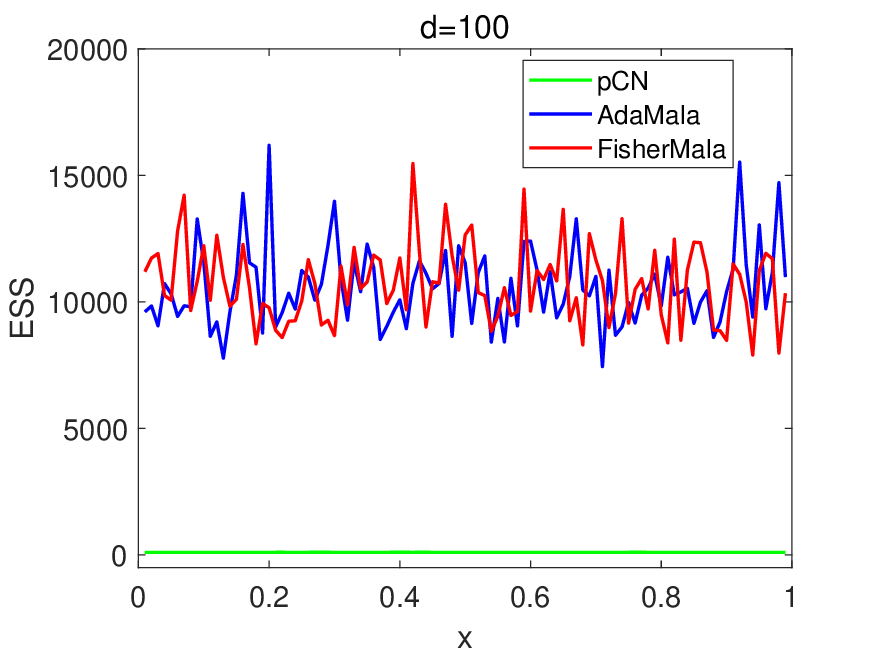}
}
\subfigure[d=600]{
\label{fig2:subfig4}
\includegraphics[width=0.43\textwidth,trim=4 0 35 18,clip]{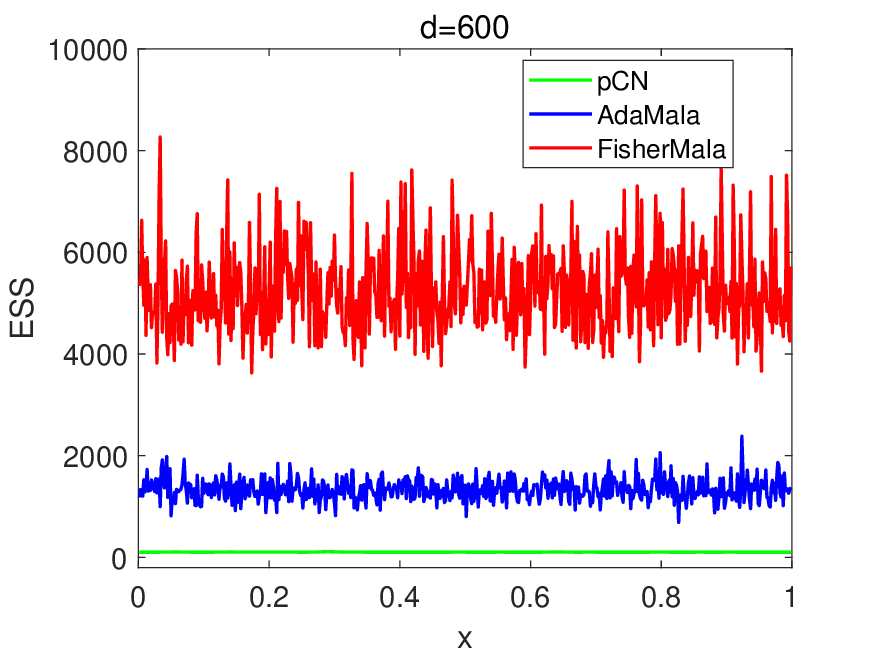}
}
\caption{Point-wise ESS for unknown function $f(x)$ with $lags=500$.
}
\label{Isess}
\end{figure}
\begin{figure}[htpb]
\centering
\subfigure[d=100]{
\label{esstime1}
\includegraphics[width=0.43\textwidth,trim=4 0 25 18,clip]{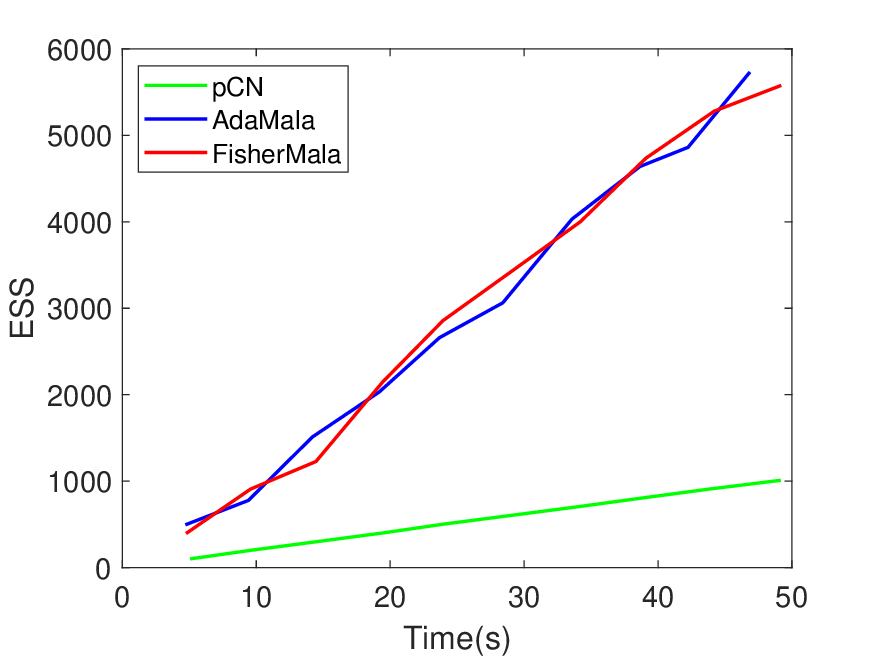}
}
\subfigure[d=600]{
\label{esstime2}
\includegraphics[width=0.43\textwidth,trim=4 0 25 18,clip]{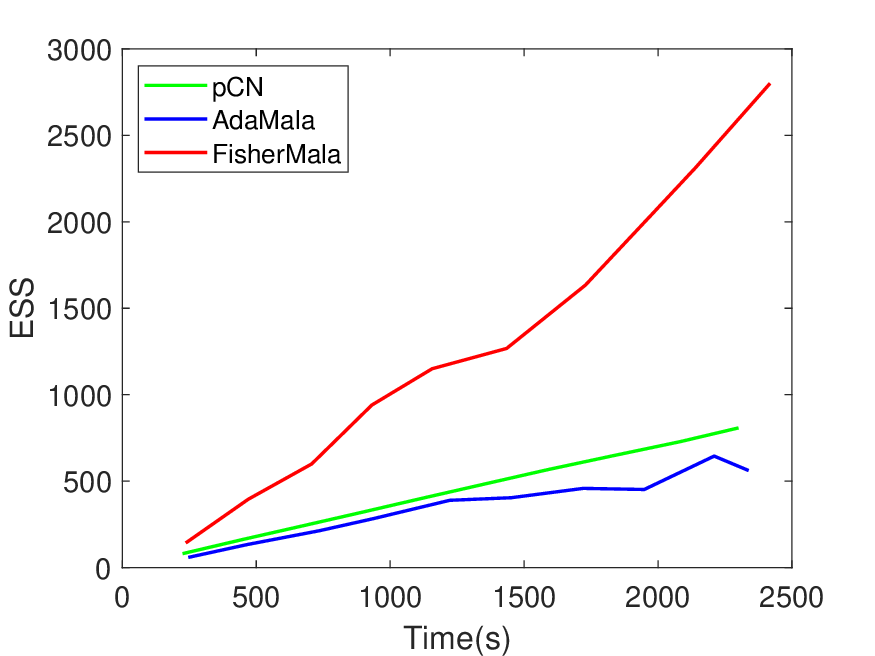}
}
\caption{The increase curves of ESS with respect to time at fixed $lags=500$.
}
\label{Isesstime}
\end{figure}

Next, we will further compare the sampling efficiency of the different algorithms. The autocorrelation function (ACF) and the effective sample size (ESS) are widely used criteria to assess the efficiency of sampler (details on the computation of the ACF and ESS can be found in Appendix \ref{ESSACF}). Utilizing the samples in collection phase after burn-in phase, we calculate the ACF and ESS for all algorithms.
Figure \ref{Isacf} displays the ACF curves for the state space point $x_p=0.5$ of the unknown function $f(x)$ with $lag=500$. Point-wise ESS estimates for the source term $f(x)$ are displayed in Figure \ref{Isess}.
For $d=100$, both FisherMala and AdaMala achieve similar ESS values, significantly outperforming pCN. However, for $d=600$, the ESS values of FisherMala are several orders of magnitude higher than the other two algorithms uniformly well in each dimensions.
This is because the samples from the pCN algorithm and AdaMala have higher correlation, as shown in Figure \ref{Isacf}.
To further compare efficiency, we compute the curves of ESS with respect to time (second) and plot the results in Figure \ref{Isesstime}. We can see that the FisherMala achieves a lager ESS  in the least amount of time, particularly in the $d=600$ case. Although pCN algorithm is gradient-free and requires minimal computation time per sample, its high autocorrelation brings significantly longer time to achieve the same ESS values. For $d=100$, AdaMala performs similarly to FisherMala. While, for $d=600$, FisherMALA has the highest ESS/time efficiency.
\begin{figure}[htpb]
\centering
{
\label{Iserr:fig1}
\includegraphics[width=0.65\textwidth,trim=3 3 25 17,clip]{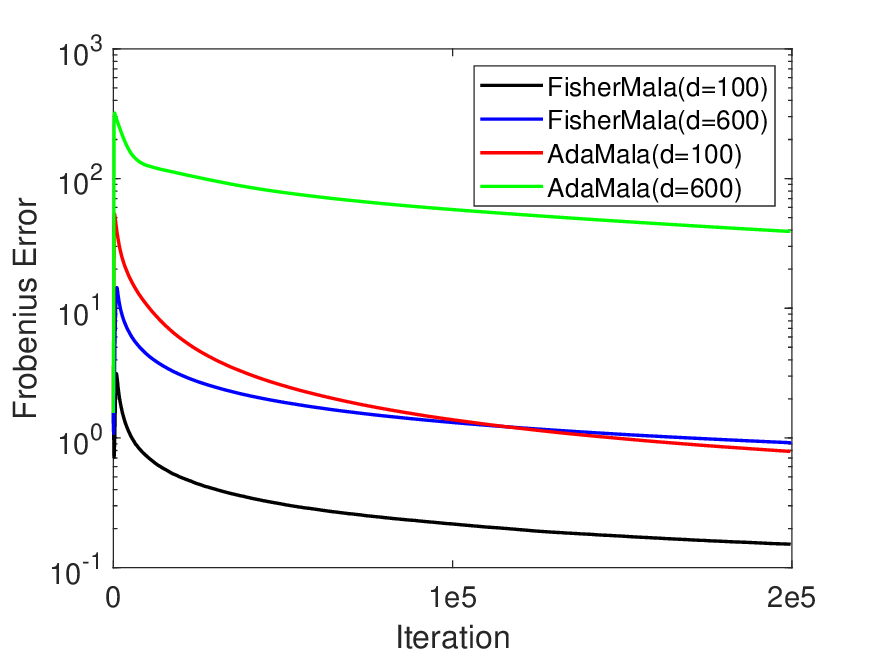}
}
\caption{The Frobenius norm across iterations with different adaptive preconditioner.}
\label{Iserr}
\end{figure}
Within the linear Bayesian inverse problem framework, if the prior is Gaussian distribution, and the posterior distribution is also a Gaussian distribution with the covariance matrix given by equation \eqref{cpost} as previous discussion in Section \ref{baye}. Obviously, for a multivariate Gaussian $\pi(x)=N(\mu_{post}, C_{post})$, it holds the inverse Fisher matrix equals the covariance matrix, i.e.  $\mathcal{I}^{-1}=C_{post}$. We adopt the same prior distribution $\pi_0=N(0,1.5 I_d)$ for both FisherMala and AdaMala, and compare the convergence speed of the preconditioner to the posterior covariance matrix.
To quantify convergence, we compute the Frobenius norm of the normalized matrix $\|\widetilde{M}_n - \widetilde{C}_{post}\|$ across adaptation iterations $n$, where $\widetilde{B}$ denotes
 trace-normalized matrices, i.e.\  $\widetilde{B} = B / (\frac{\text{trace}(B)}{d})$,  and $\widetilde{C}_{post}$ is the normalized posterior covariance matrix. Figure \ref{Iserr} shows the error decay curves for dimensions $d = 100$ and $d = 600$. For better visibility, the $y$ axis is shown in log scale.
For $ d = 100$, both algorithms exhibit rapid convergence. However, for $ d = 600$, AdaMala converges significantly more slowly. The transient error increase during initial iterations due to insufficient exploration of the posterior distribution. Afterward, the errors decay rapidly and become stable. This phenomenon is consistent with the theoretical result in proposition \ref{prop:convergence}. The FisherMala converges faster than AdaMala for both lower and higher dimensional settings.

From all the results, we conclude that FisherMala outperforms the other algorithms, particularly in high-dimensional settings. Its superior ESS, lower autocorrelation, and faster convergence make it an efficient and robust sampler for solving Bayesian inverse problems.

\subsection{The parameter identification problem}\label{example2}
The parameter identification problem represents a class of nonlinear inverse problems that are complex than the linear problem discussed in example \ref{example1}.
Specifically, we consider the reconstruction of parameter $q$ in the following Neumann boundary value problem
\begin{equation}\label{problem2}
\left\{
\begin{array}{ll}
-\Delta u+qu=f,  &{in\ \Omega,}\\
\dfrac{\partial u}{ \partial n}=0,  &{on\ \partial \Omega.}\\
\end{array} \right.
\end{equation}
where $q\in L^2(\Omega)$ denotes the parameter field what we aim to identify, $u\in H^1(\Omega)$ represents the state variable, and $f$ is a given source term.
Given the source term $f$ in problem \eqref{problem2}, our objective is recovering the coefficient $q(x)$ using interior measurements $u(x) |_{\Omega}$, following the framework in \cite{Gu2021}. Thus, this problem will be reduced to the following nonlinear problem
$$y=\mathcal{F}(q)+\eta,$$ where $\eta$ represents measurement noise.
In our work, we apply the perturbation difference approximation method developed in \cite{Li2008} to compute the Fr\'{e}chet derivatives of the observed operator $\mathcal{F}$.

\begin{figure}[htpb]
\centering
\subfigure[]{
\label{Idq}
\includegraphics[width=0.43\textwidth,trim=10 4 35 15,clip]{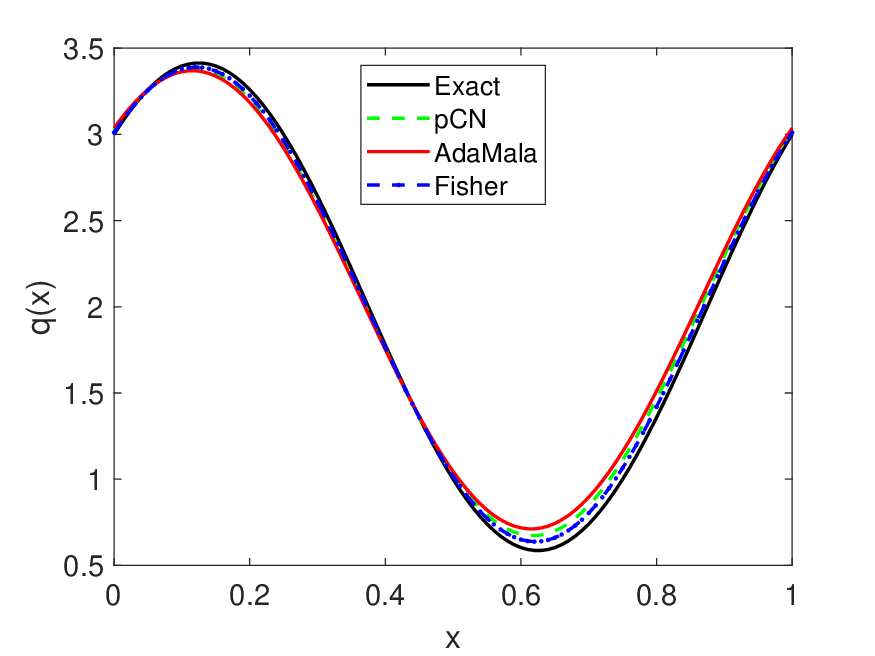}
}
\subfigure[]{
\label{Ids2}
\includegraphics[width=0.43\textwidth,trim=5 2 36 7,clip]{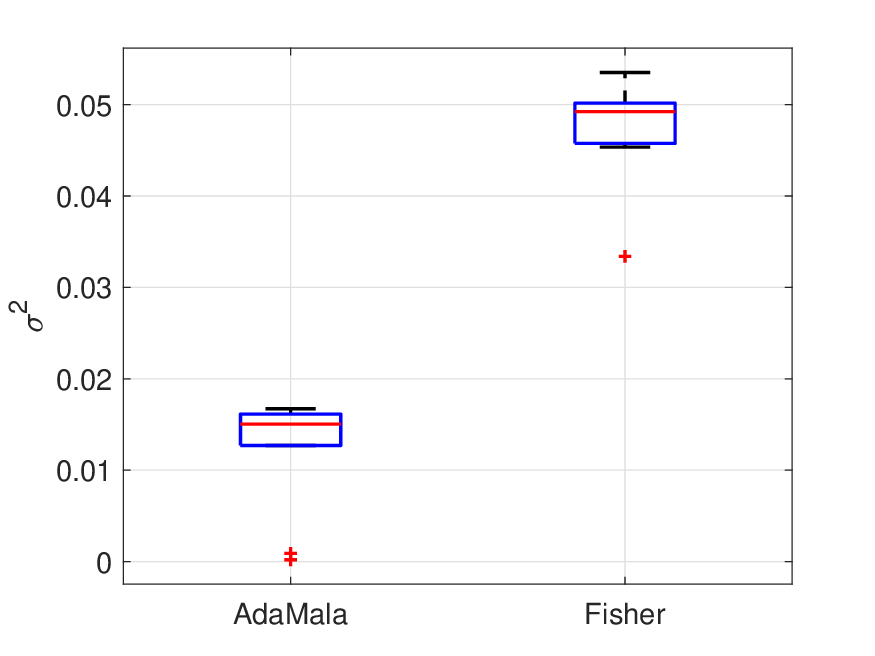}
}
\vspace{-0.3cm}
\caption{(a): The exact solution and reconstruction results for all samplers; (b): The estimated values of stepsize $\sigma^2$ of FisherMala and AdaMala algorithms.}
\label{Idqs2}
\end{figure}

\begin{figure}[htpb]
\centering
\includegraphics[width=0.9\textwidth,trim=15 0 5 5,clip]{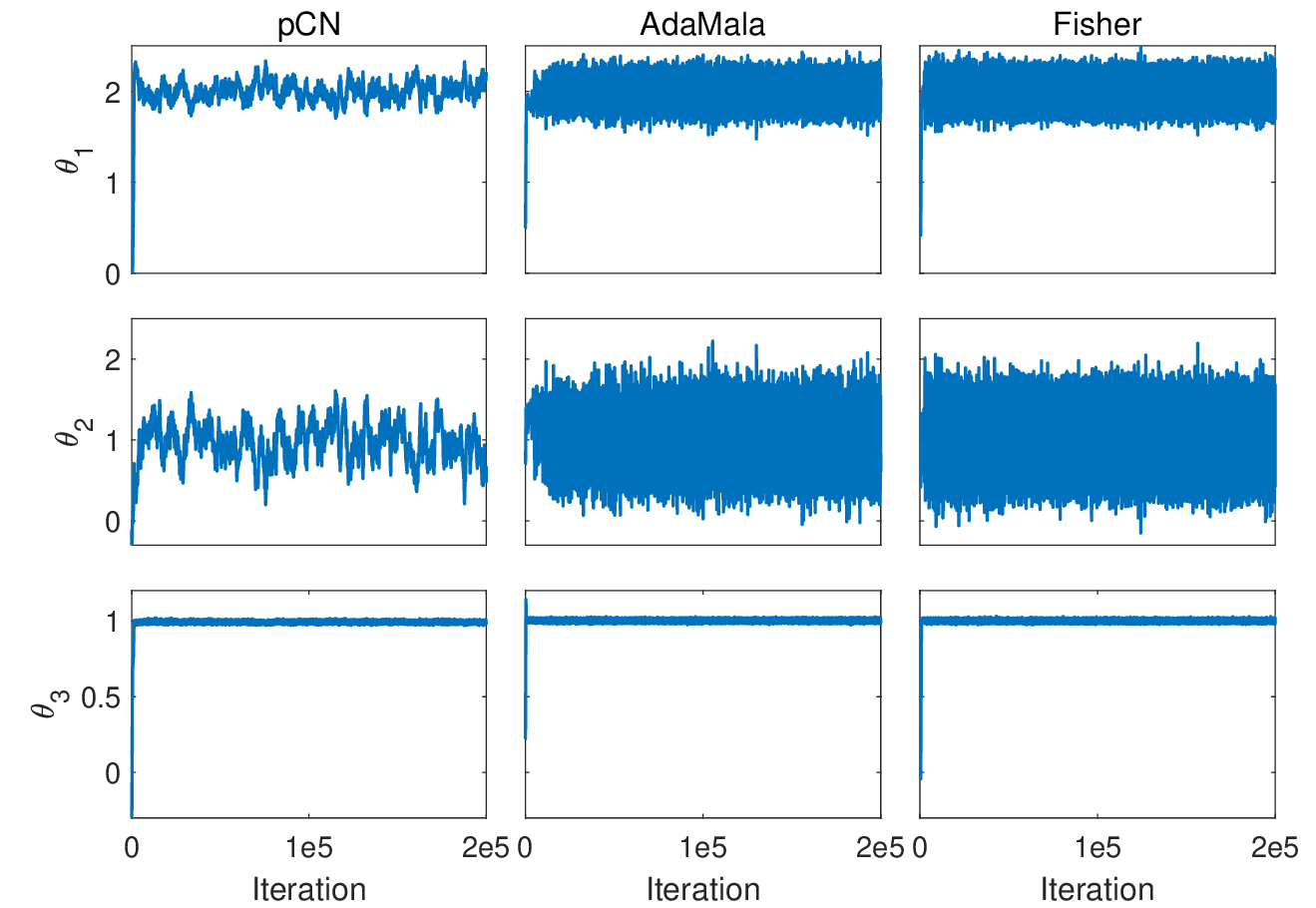}
\caption{The trace plots for all the coefficient elements of $\theta$.}
\label{Idtrace}
\end{figure}

\begin{figure}[htpb]
\centering
\subfigure[pCN]{
\label{Idpos:fig1}
\includegraphics[width=0.31\textwidth,trim=3 0 25 5,clip]{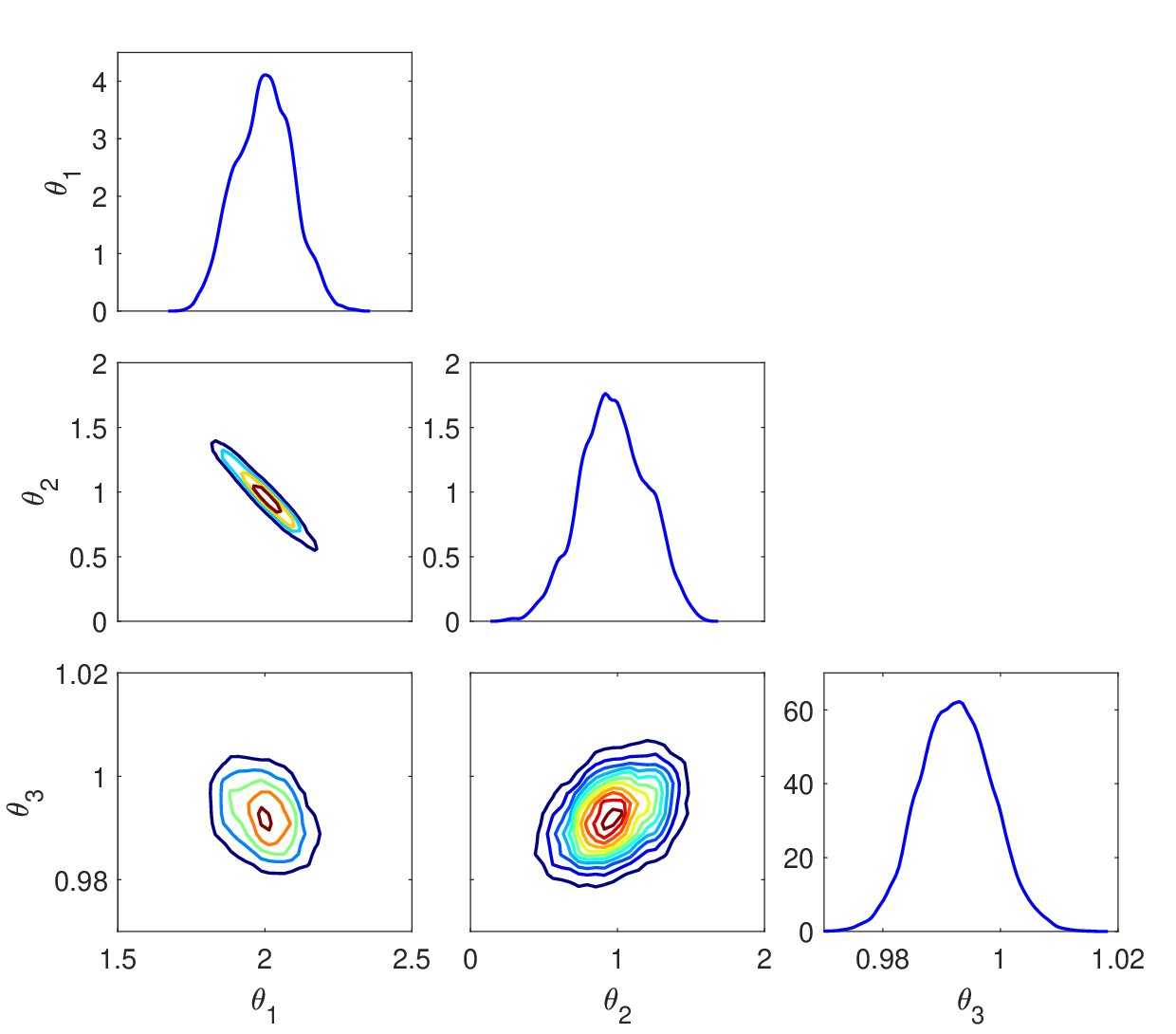}
}
\subfigure[AdaMala]{
\label{Idpos:fig2}
\includegraphics[width=0.31\textwidth,trim=3 0 25 5,clip]{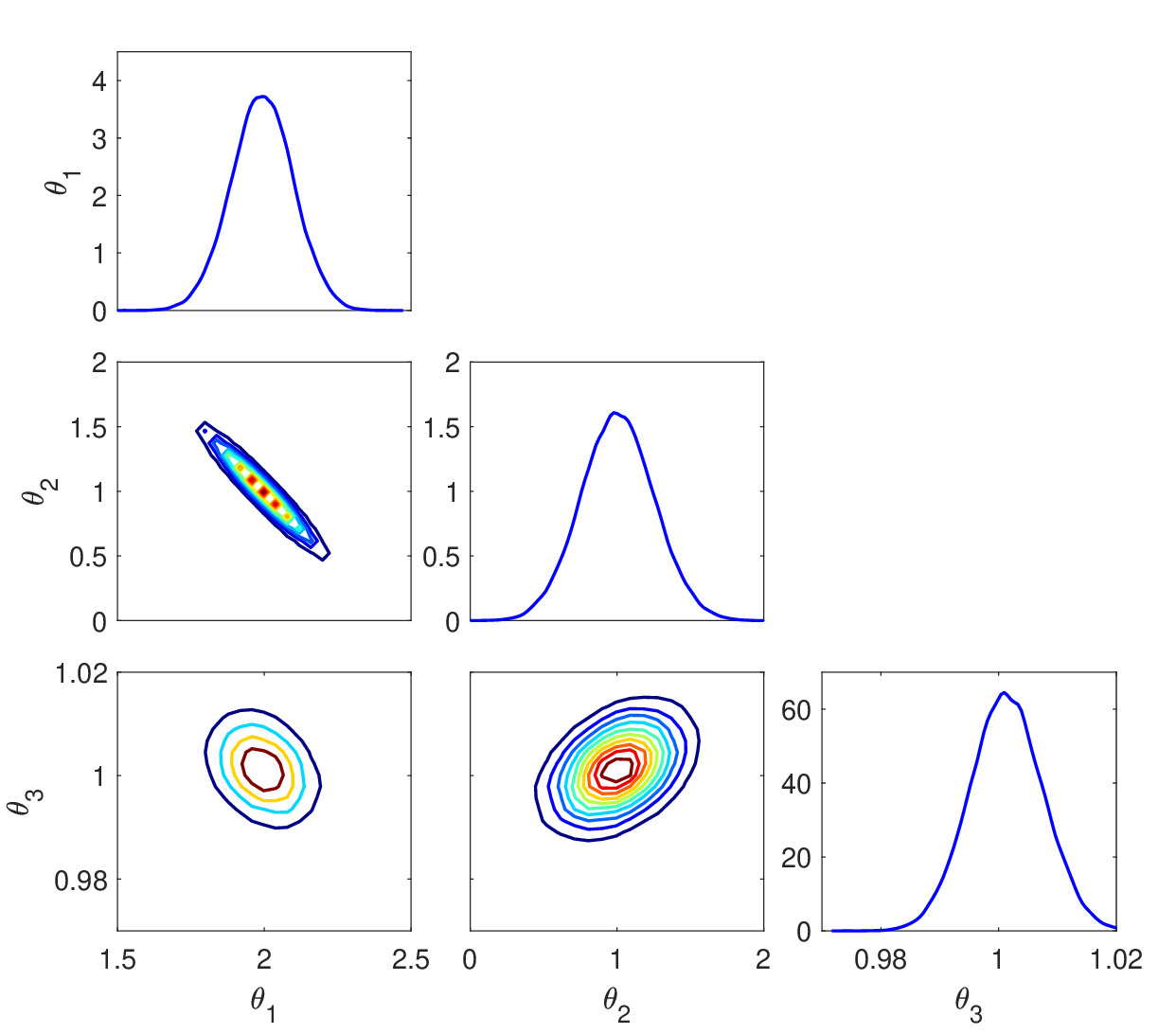}
}
\subfigure[FisherMala]{
\label{Idpos:fig3}
\includegraphics[width=0.31\textwidth,trim=3 0 25 5,clip]{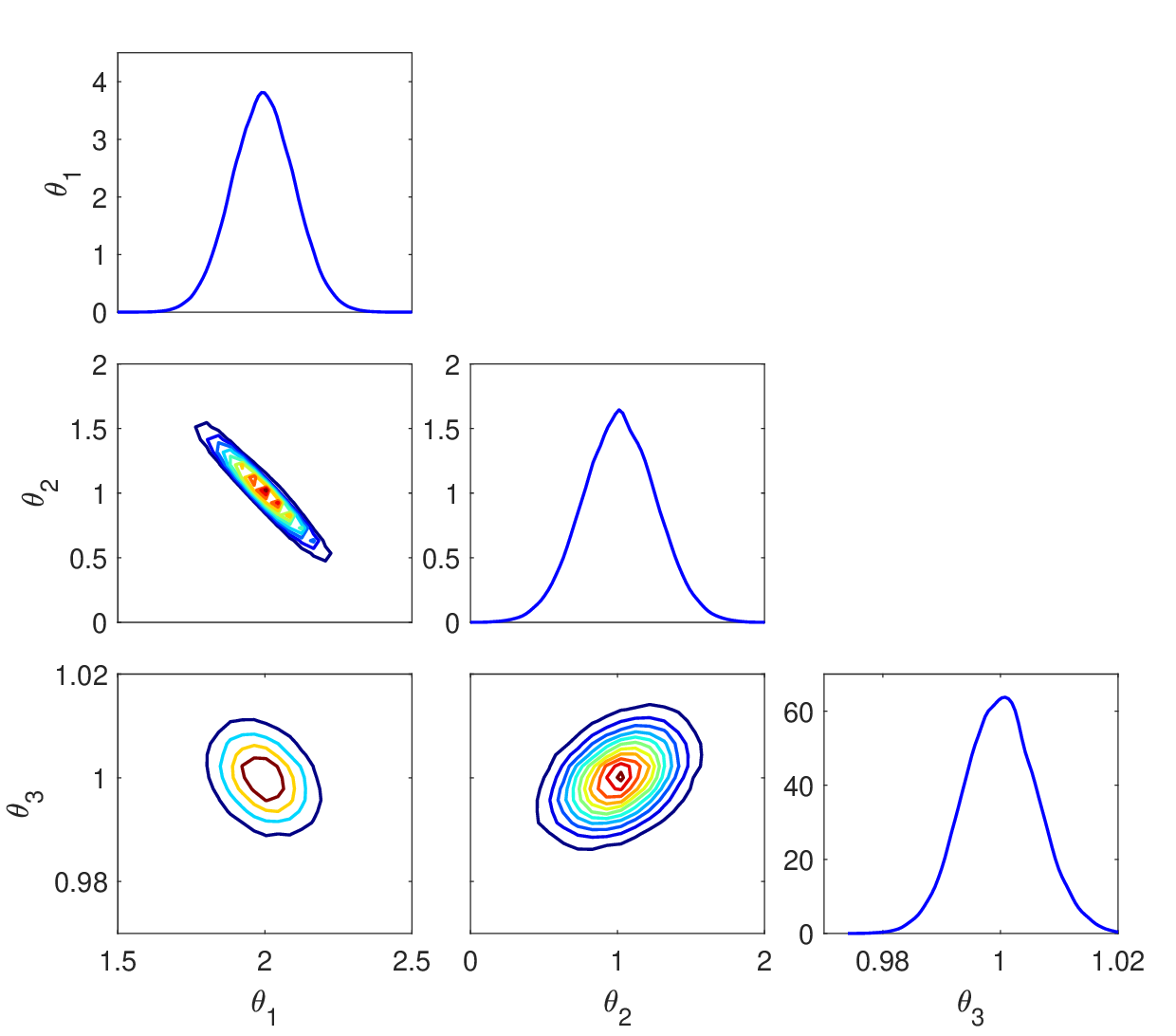}
}
\caption{The 1-D and 2-D posterior marginal densities. }
\label{Idpos}
\end{figure}

For our numerical implementations, we consider $\Omega = (0,1)$, and the source $f$ is given by
\[f(x)=q(x)\cos(\pi x)+\pi^2 \cos(\pi x).\]
The exact parameter field $q(x)$ is expressed in the trigonometric basis $$\{1, \sin(2\pi x), \cos(2\pi x), \cdots, \sin(2N\pi x), \cos(2N\pi x)\}$$ of $L^2(\Omega)$ as following
\begin{equation}\label{q}
q(x)= 2+\sin(2\pi x)+\cos(2\pi x),
\end{equation}
where $\theta = (\theta_1, \theta_2, \theta_3) = (2, 1, 1)$ represents the target coefficients as \cite{Li2008}. The differential equation \eqref{problem2} is discretized using a uniform grid with $N_x=100$ and approximately solve with second order centered difference scheme. For all experiments, we take the prior distribution as $\pi_0=N(0, 0.1 I_d)$.
We draw $2\times10^5$ samples from posterior distribution, discarding the first $1\times10^5$ as burn-in process.
With a fixed error level $\epsilon=0.01$, we compare the performances of different algorithms. The inversion results and performances for the three methods are summarized in Table \ref{qtab}, while Figures \ref{Idqs2}, \ref{Idtrace}, and \ref{Idpos} illustrate the findings.
\begin{table}[htpb]
\centering
\caption{The performance in Example 2.}
\setlength{\tabcolsep}{5pt} 
\renewcommand{\arraystretch}{1.3}
\begin{tabular}{cccc}
\toprule
{} &{Inversion results} & {err(\%)} & {ESS}\\
\midrule
pCN & (2.026, 0.914, 0.998) & 3.66 & (128, 124, 1265) \\
AdaMala &(2.040, 0.880, 0.995)  & 5.17 & (16241, 16035, 36374) \\
FisherMala & (2.014, 0.948, 0.996)  & 2.17  & (\textbf{57246, 53032, 56561}) \\
\bottomrule
\end{tabular}
\label{qtab}
\end{table}
The inversion results presented in table \ref{qtab} and Figure \ref{Idq} indicate that FisherMala achieves superior reconstruction accuracy. Fixed the $lag=500$, we compute the ESS use the samples after burn-in phase. The results in table \ref{qtab} reveal that FisherMala has a higher ESS, which can be orders of magnitude larger than pCN algorithms and significantly larger than AdaMala. The step size for the pCN algorithm is set to $\beta=0.055$, yielding an acceptance rate approximate $25\%$. In this numerical experiment, it is crucial to ensure that the step size $\sigma^2$ exceeds $1\times10^{-7}$ for both FisherMala and AdaMala. Through multiple experiments, we obtain the stepsize parameter $\sigma^2$ close to $0.05$ and $0.015$ for FisherMala and AdaMala, respectively, as shown in Figure \ref{Ids2}. 
Figure \ref{Idtrace} displays the sample pathes for each element of coefficient $\theta$. We can see that the pCN algorithm demonstrates the poorest mixing times, while FisherMala performs best for the coefficients $\theta$. Consequently, the FisherMala efficiently explores the state space.
%
We also visualize the 1-dimensional (1-D) and 2-dimensional (2-D) marginal posterior probability densities for the parameters $\theta$ across the three methods in Figure \ref{Idpos}. The shape of the 2-D marginal densities indicate apparent correlations among certain parameters, particularly between $\theta_1$ and $\theta_2$. These also show that the posterior distributions are  non-Gaussian distribution and verify the complexity of the nonlinear inverse problems.

\begin{figure}[htpb]
\centering
\subfigure[pCN]{
\label{Idnsy:fig1}
\includegraphics[width=0.31\textwidth,trim=7 0 32 5,clip]{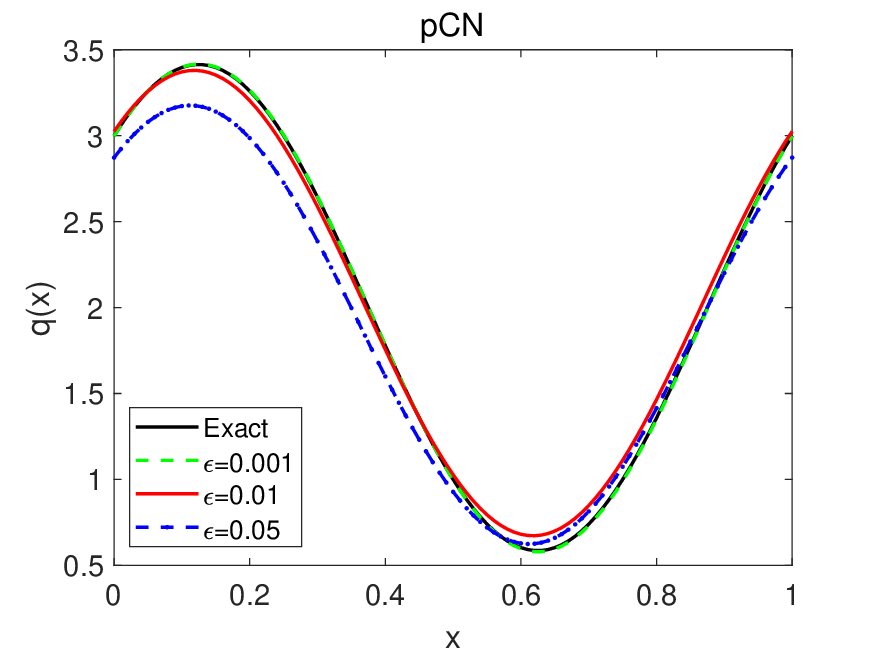}
}
\subfigure[AdaMala]{
\label{IIdnsy:fig2}
\includegraphics[width=0.31\textwidth,trim=7 0 32 5,clip]{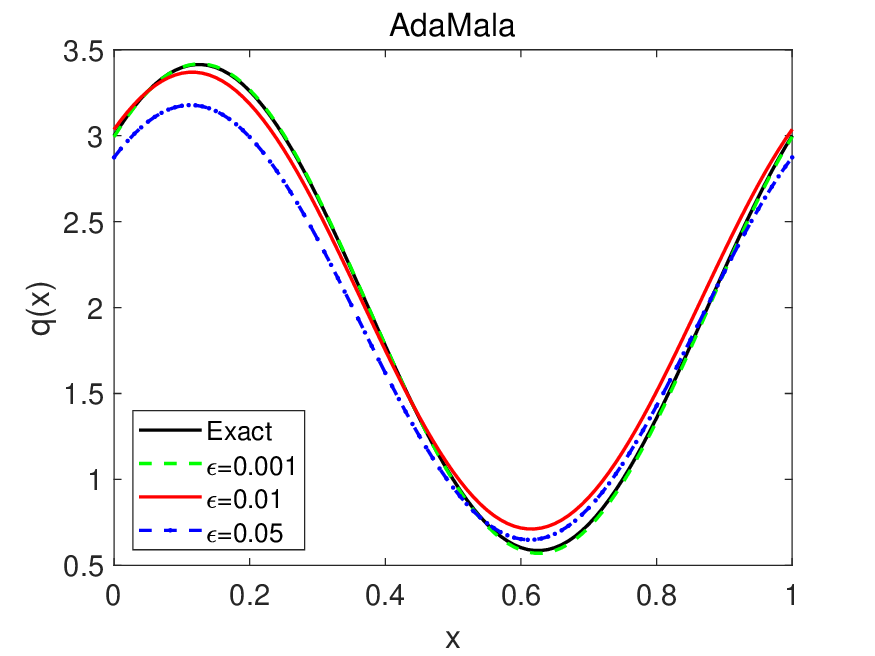}
}
\subfigure[FisherMala]{
\label{Idnsy:fig3}
\includegraphics[width=0.31\textwidth,trim=7 0 32 5,clip]{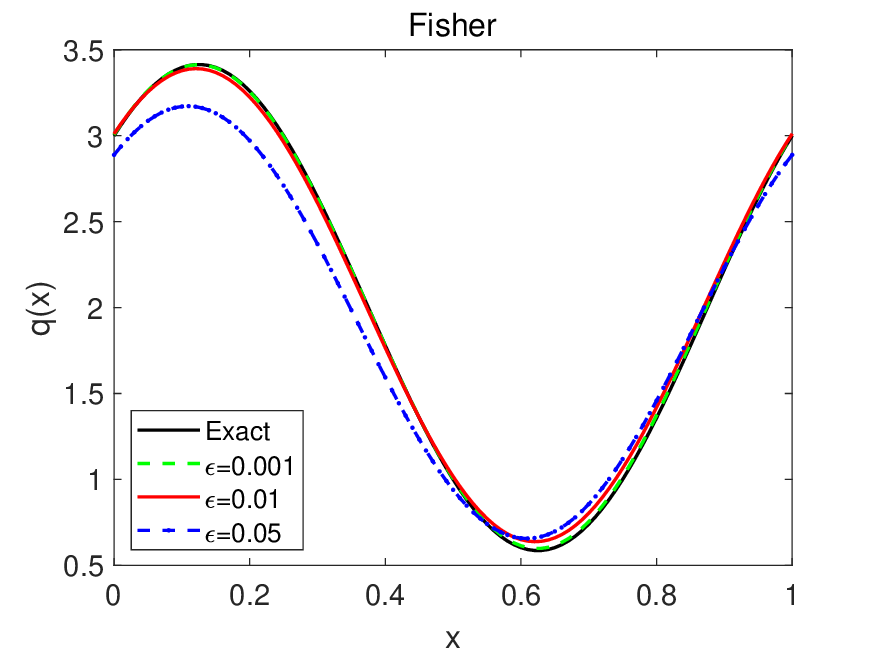}
}
\caption{The reconstructions with different noisy levels. }
\label{Idnsy}
\end{figure}

Finally, we assess the stability of the methods against various noise levels of $0.001$, $0.01$ and $0.05$. The results for the three samplers are displayed in Figure \ref{Idnsy}. As the noise level decreasing, the inversion results converge towards the exact results for all the methods. However, FisherMala algorithm is slightly more accurate. Overall, we can conclude that FisherMala is better than all other algorithms from these findings.

\subsection{An advection-diffusion-reaction (ADR) inverse problem}\label{example3}
\begin{figure}[htpb]
\centering
{
\label{ADR_m0}
\includegraphics[width=0.23\textwidth,trim=10 5 20 5,clip]{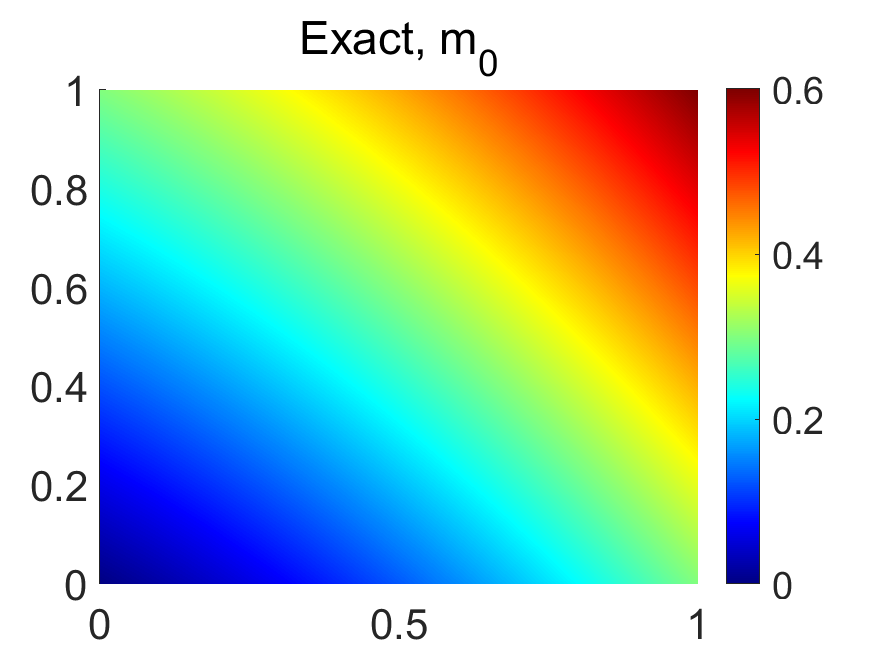}
}
{
\label{ADR_pCNm0}
\includegraphics[width=0.23\textwidth,trim=10 5 20 5,clip]{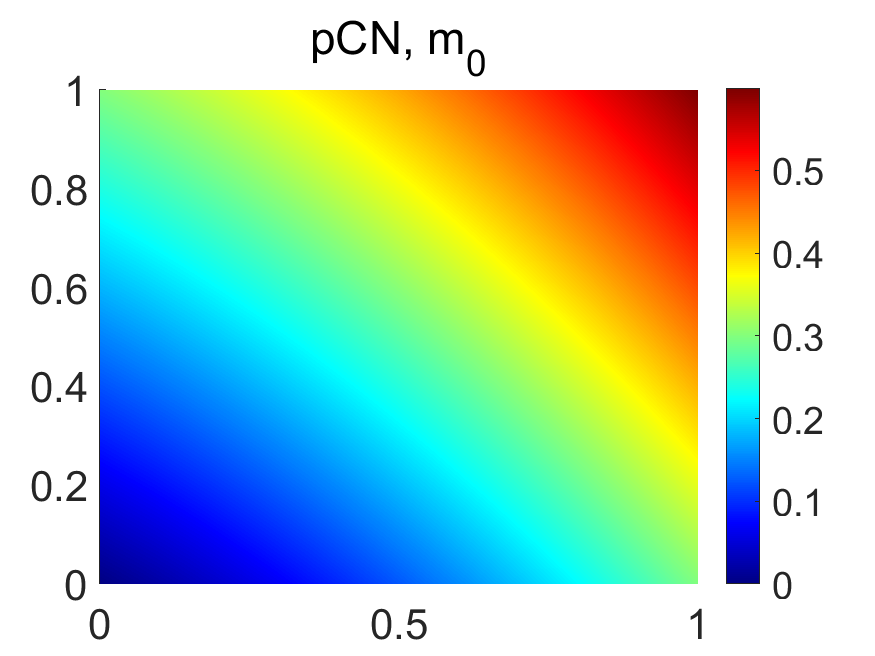}
}
{
\label{ADR_malam0}
\includegraphics[width=0.23\textwidth,trim=10 0 20 5,clip]{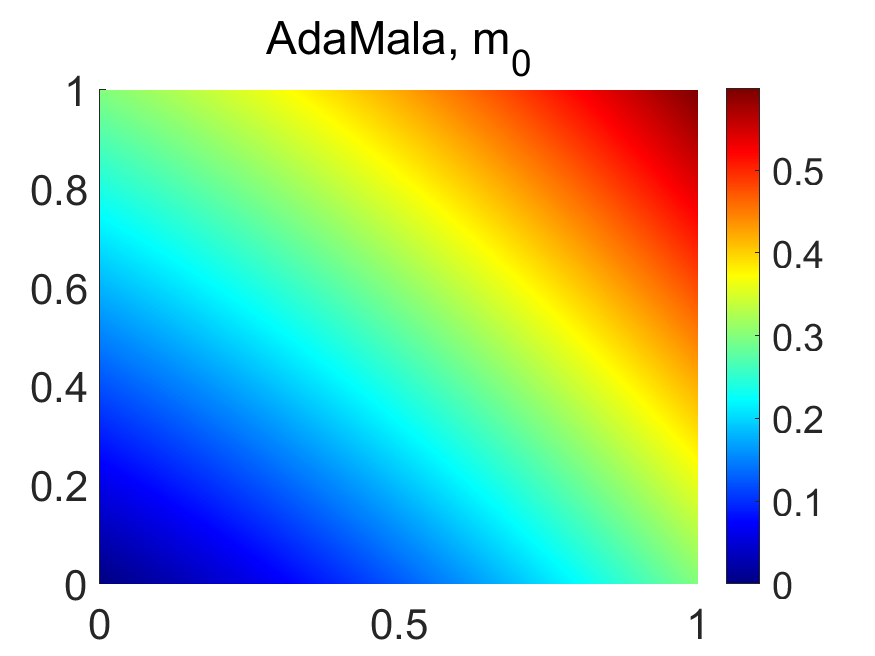}
}
{
\label{ADR_fisherm0}
\includegraphics[width=0.23\textwidth,trim=10 0 20 5,clip]{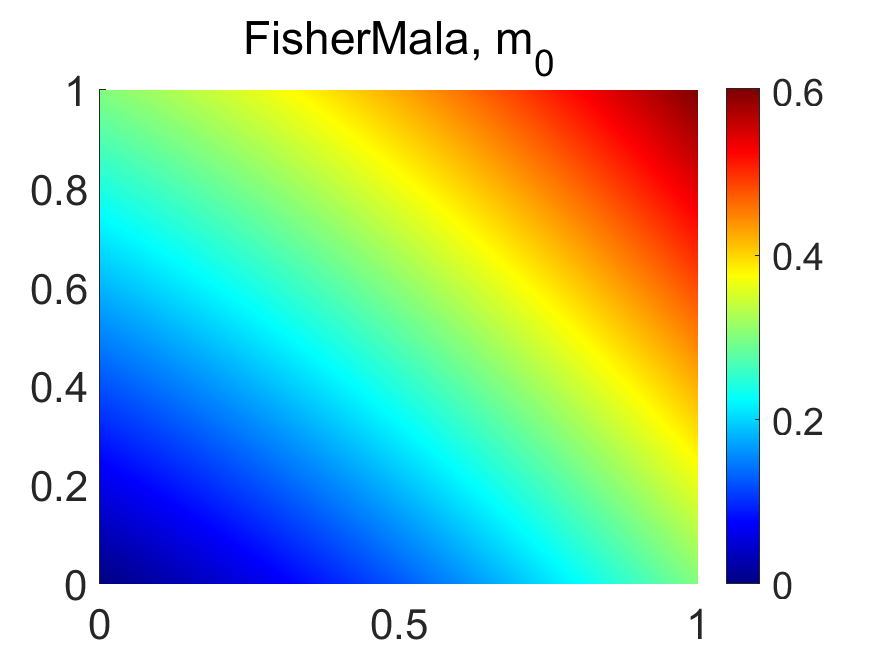}
}
\vspace{0.5cm}\\
{
\label{ADR_mk}
\includegraphics[width=0.23\textwidth,trim=10 0 5 5,clip]{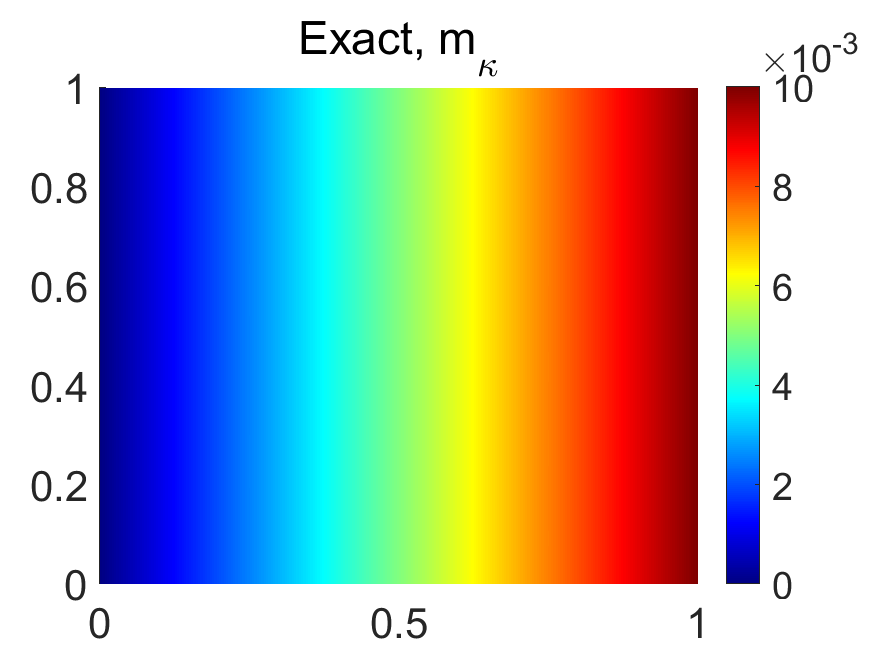}
}
{
\label{ADR_pCNmk}
\includegraphics[width=0.23\textwidth,trim=10 5 5 5,clip]{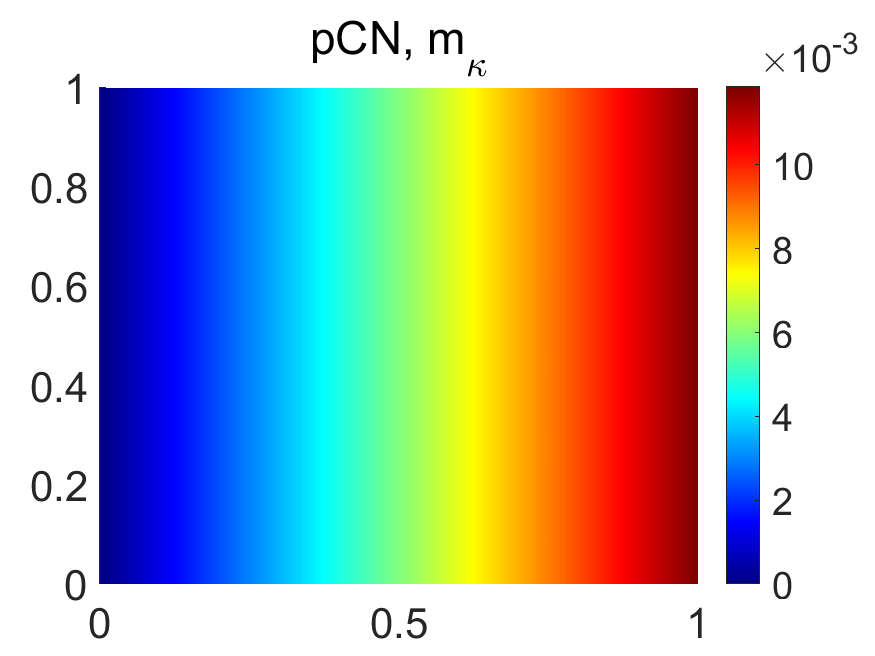}
}
{
\label{ADR_malamk}
\includegraphics[width=0.23\textwidth,trim=10 5 5 5,clip]{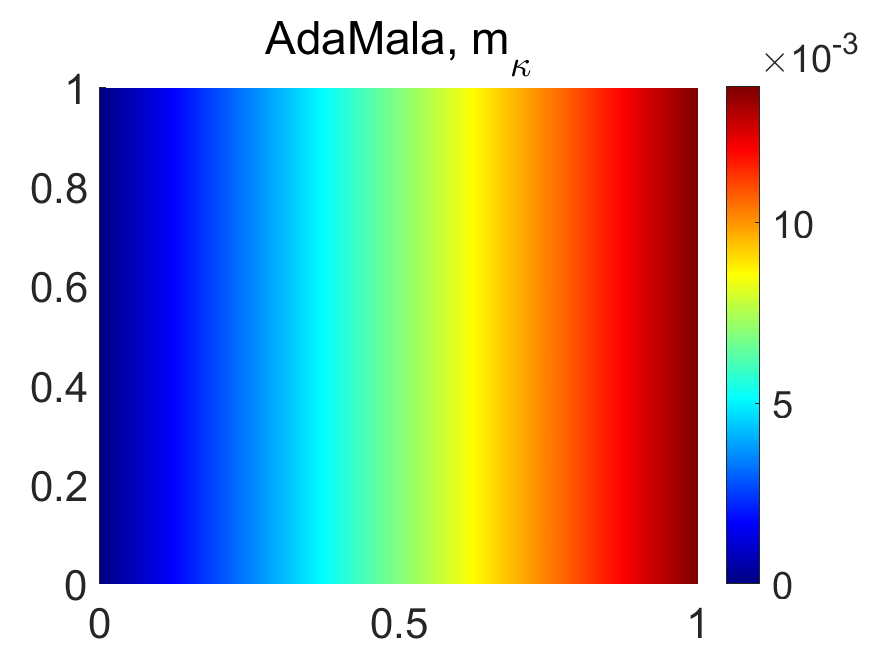}
}
{
\label{ADR_fishermk}
\includegraphics[width=0.23\textwidth,trim=10 5 5 5,clip]{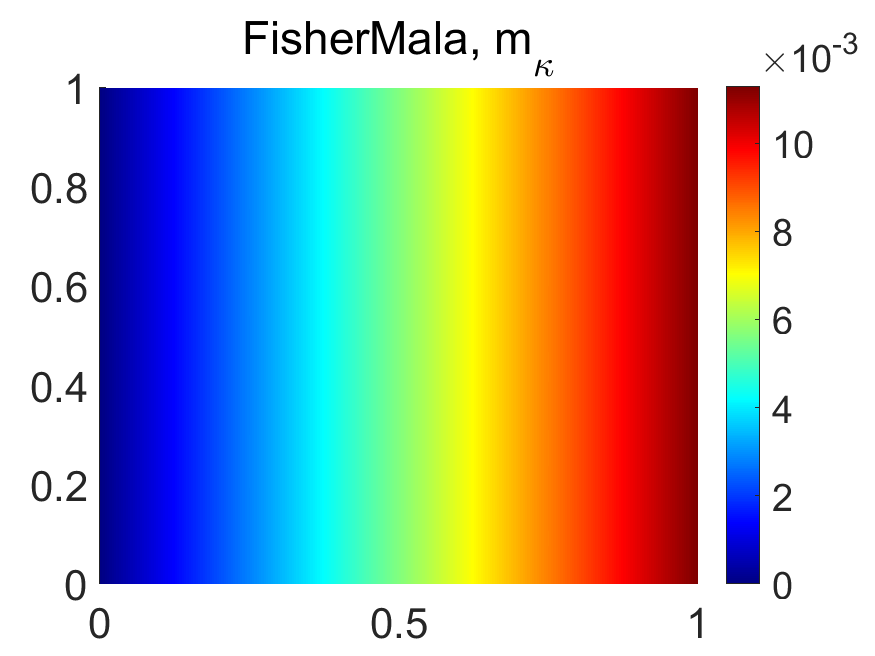}
}
\caption{The exact and reconstruction results.
}
\label{ADR}
\end{figure}

\begin{figure}[htpb]
  \center
 {
\label{ADR_acf1}
\includegraphics[width=0.43\textwidth,trim=5 0 25 7,clip]{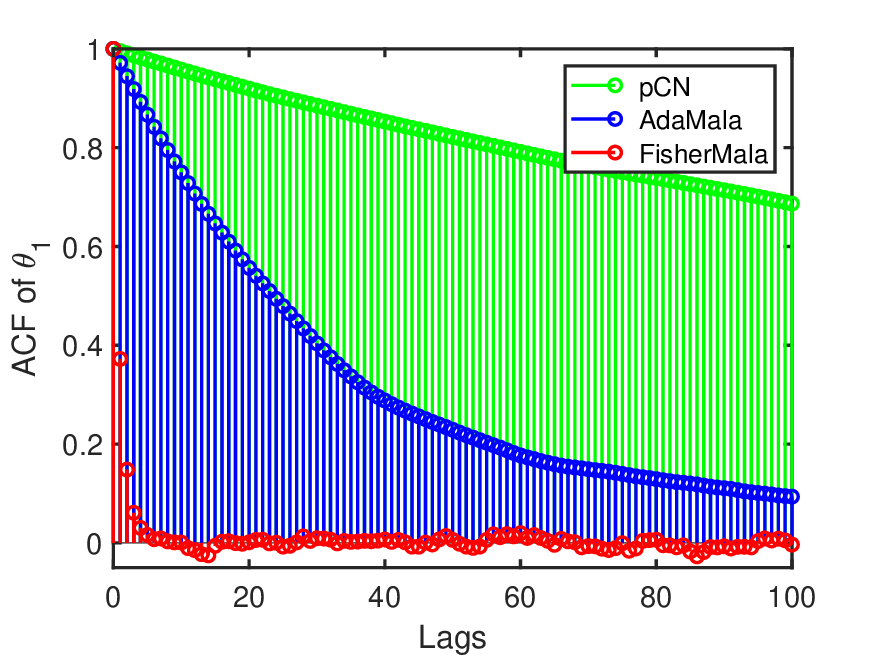}
}
{
\label{ADR_acf2}
\includegraphics[width=0.43\textwidth,trim=5 0 25 7,clip]{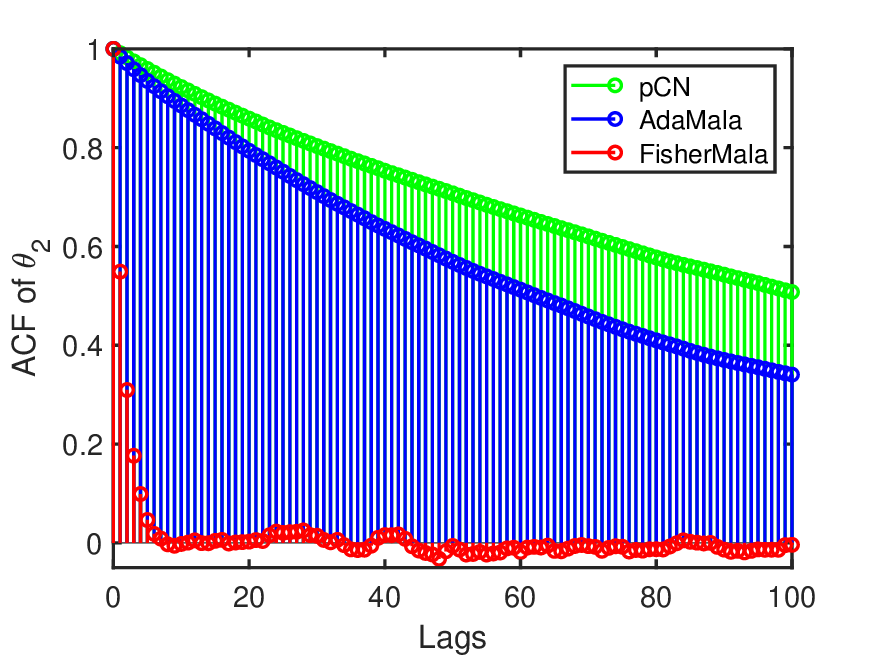}
}
{
\label{ADR_acf3}
\includegraphics[width=0.43\textwidth,trim=5 0 25 7,clip]{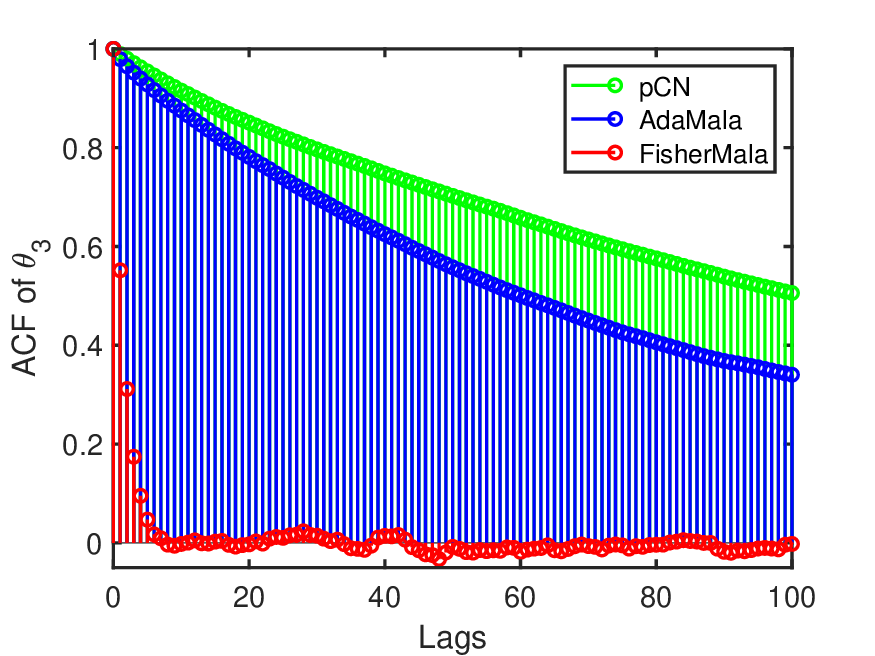}
}
{
\label{ADR_acf4}
\includegraphics[width=0.43\textwidth,trim=5 0 25 7,clip]{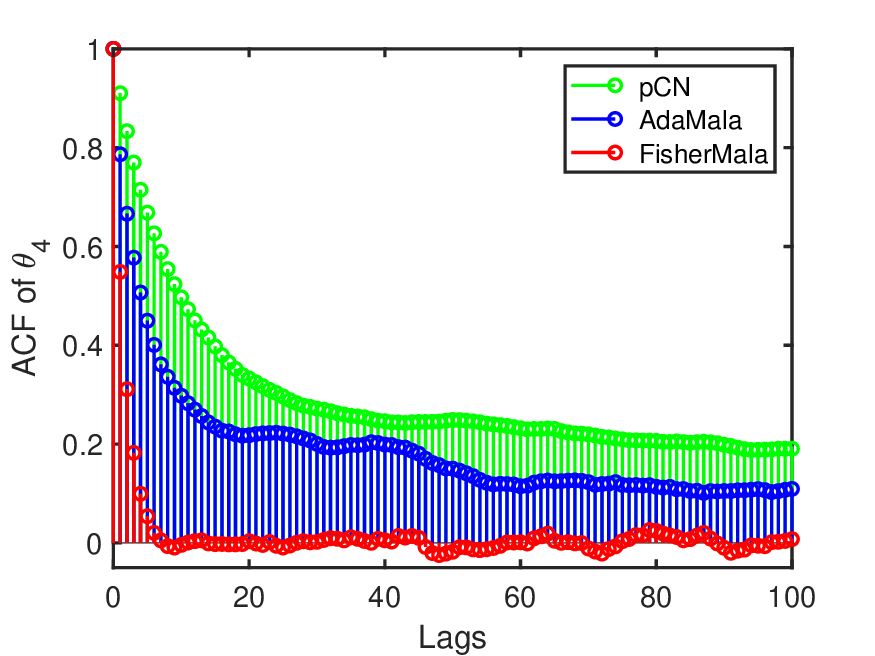}
}
  \caption{The autocorrelation functions (ACF) plot with $lags=100$.}
  \label{ADRacf}
  \end{figure}

\begin{figure}[htpb]
\centering
\subfigure[]
{
\label{ADRess}
\includegraphics[width=0.43\textwidth,trim=5 4 35 15,clip]{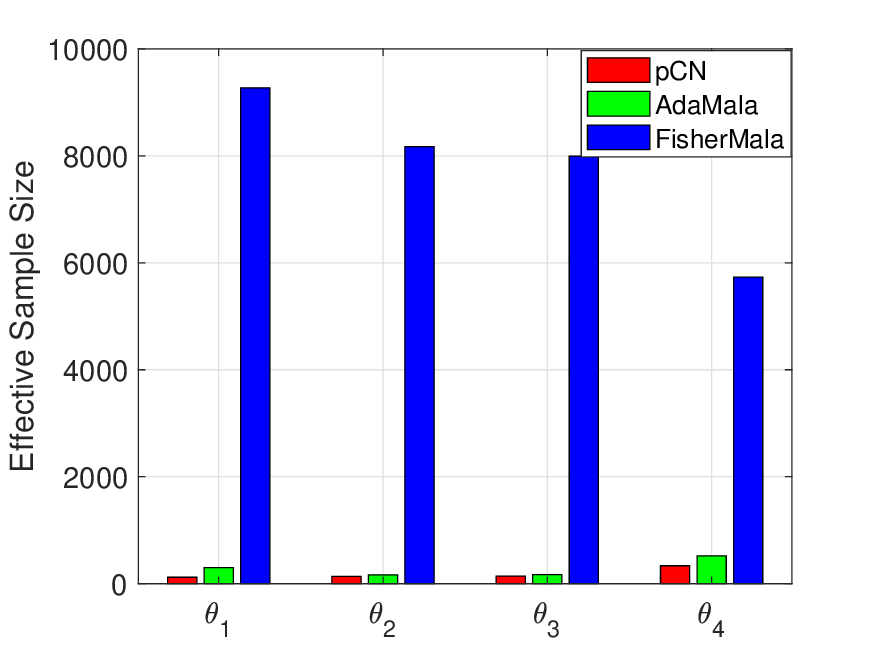}
}
\subfigure[]
{
\label{ADRCI}
\includegraphics[width=0.43\textwidth,trim=5 2 36 7,clip]{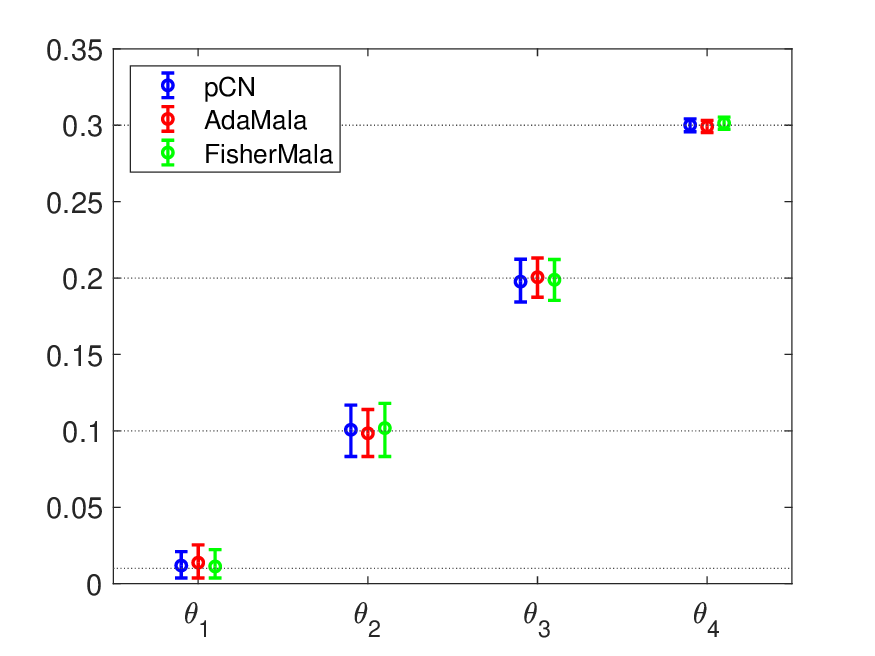}
}
\vspace{-0.3cm}
\caption{(a): The effective sample size (ESS); (b): The $95\%$ confidence intervals (CI) for different methods, the dotted line represents the exact parameters, colored 'o' represents mean values.}
\label{ADRess-CI}
\end{figure}

The simultaneous reconstruction of the diffusion coefficient field $m_\kappa$ and initial condition $m_0$ problem is significant challenges due to nonlinear, as discussed in \cite{Ghattas2021}.
In this example, we consider a more complex inverse problem governed by the two-dimensional nonlinear time-dependent ADR equation:
\begin{equation}\label{problem3}
\left\{
\begin{array}{ll}
\dfrac{\partial u}{\partial t} +\mathbf{v}\cdot \nabla u-\nabla\cdot (m_\kappa \nabla u)+cu^3=f,  &{in\ \Omega\times(0,T),}\vspace{0.1cm}\\
m_{\kappa}\dfrac{\partial u}{ \partial n}=0,  &{on\ \partial \Omega\times(0,T),}\vspace{0.14cm}\\
u|_{t=0}=m_0, &{in\ \Omega,}
\end{array} \right.
\end{equation}
where $\mathbf{v}$, $m_\kappa$ and $c$ denote the velocity field, diffusion coefficient and reaction coefficient, respectively. Set the computational domain $\Omega=(0,1)\times(0,1)$, the time-dependent advection velocity field $\mathbf{v(t)}=(\cos(t),\sin(t))$ and the reaction coefficient $c=1$. Synthetic measurements are generated from the interior state field $u(x,y,T)|_{\Omega}$ at terminal time $T=0.05$.
Source term $f(x,y)$ is given by $$f(x,y)=\exp\{((x-0.5)^2+(y-0.5)^2)/0.9^2\}$$.
The exact solution of the diffusion coefficient $m_\kappa$ and the initial condition $m_0$ are defined by polynomial basis $\{x,y,x^2,\cdots, xy^{N-1},y^{N}\}$ as following
\begin{align*}
&m_\kappa(x,y)=0.01x,\\
&m_0(x,y)=0.1x^2+0.2x+0.3y,
\end{align*}
and we aim to recover the coefficients $\theta=(\theta_1,\theta_2,\theta_3,\theta_4)=(0.01,0.1,0.2,0.3)$.
The nonlinear differential equation \eqref{problem3}
is solved by the finite element method (FEM) with the codes from \cite{Quarteroni2015}, and Newton Method to linearization the nonlinear term $u^3$. We choose the linear basic function and the maximal space elements size $0.1$. Discretize the time domain with uniform grids $T/20$. Fixed the error level $\epsilon=0.01$, the prior distribution is $\pi_0=N(0,C)$ with $C=[0.01~ 0~ 0~ 0;~0~ 0.1~ 0~ 0;~0~ 0~ 0.1 ~0;~0~ 0~ 0 ~0.1]$. For all experiments and samplers, we consider $2\times 10^4$ burn-in iterations and $2\times 10^4$ iterations for collecting samples. The initial samples are from the prior distribution, and the first component of the initial points $\theta_0$ must be positive. Then, we show the numerical simulation results and compare the performance of different samplers in Figure \ref{ADR}-\ref{ADRess-CI}. Here, ESS and ACF are computed from the samples after burn-in process with fixed $lag=100$.


Our results demonstrate that all sampling methods successfully approximate the target parameters, confirming the validity of the inversion framework. However, comparative analysis through Figures \ref{ADRacf} and \ref{ADRess},  we observe distinct performance characteristics. The FisherMala exhibits significantly lower autocorrelation (ACF) and higher effective sample size (ESS) values compared to both pCN and AdaMala samplers. This indicates the superior sampling efficiency of FisherMala in parameter space exploration for inverse problems, same as previous examples.
We also present further quantification of uncertainty through the $95\%$ credible interval (CI) for the reconstruction parameters in Figure \ref{ADRCI}. The statistical analysis demonstrates that all computed confidence intervals successfully encompass the true parameter values and further highlighting the advantages of the FisherMala algorithm in sampling efficiency.

\section{Conclusion}\label{sec4}
In this research, we have presented the Fisher adaptive Metropolis adjusted Langevin algorithm (FisherMala). This algorithm utilizes the inverse Fisher information matrix as an optimal preconditioning for Langevin diffusion, by minimizing the expected squared jumped distance and adaptively learn the Fisher matrix from the history of samples. We have supplemented the convergence analysis of the online learning scheme within the framework of stochastic approximation and provided a proof of convergence rate under certain assumptions.
Furthermore, we have applied FisherMala to Bayesian inverse problems and compared its performance with the preconditioned Crank-Nicolson (pCN) and standard adaptive Metropolis adjusted Langevin algorithm (AdaMALA). Through numerical examples, we have demonstrated the superior efficiency of FisherMala over these two methods.
For nonlinear inverse problems, we have employed the perturbation difference approximation method to approximate the derivative information of the log-likelihood function. However, this method can be computationally expensive for high-dimensional problems. To address this issue, a potential direction for future work is to utilize the adjoint method to calculate the derivatives.
In this study, our focus has been on Bayesian inference with smooth unknown functions and Gaussian priors. However, for future research, we aim to extend FisherMALA by incorporating more complex priors such as the total variation (TV) prior \cite{Vogel2002}, which suited for the known with sharp jumps and discontinuities, or fractional total variation (FTV) prior \cite{wang2021,wang2023}, which suited for the known with sharp corners and combination of discontinuities and smooth regions. With this extension, FisherMALA would be able to handle a wider range of inverse problems more effectively.

\section{Acknowledgements}\label{sec5}
The work of L. Wang is supported by the Postgraduate Scientific Research Innovation Project of Hunan Province (CX20230406). The work of G. Zheng were supported by the NSF of China (12271151).

\section*{Appendix}
\appendix
\section{The initialization phase}\label{initial}
Adaptation of the parameter $\sigma^2$, occurs only during burn-in and keeps fixed at collection of samples stage.
In the burn-in process of FisherMala and AdaMala, the first $n_0=500$ iterations are used as the initialization phase, where samples are generated by standard MALA (i.e. the preconditioner $M$ is identity matrix),
which adapts only a step size $\sigma^2$. Besides, for AdaMALA, we do an additional set of $n_0=500$ iterations where standard MALA still runs and collects samples which are used to sequentially update the empirical covariance matrix $\mathcal{C}_n$ in equation \eqref{eq:adaMALA}. The purpose of this second phase is to play the role of  "warm-up" and provide a reasonable initialization for $\mathcal{C}_n$ at start in burn-in iteration as \cite{Titsias2024}.

\section{Autocorrelation function (ACF) and Effective sample size (ESS)}\label{ESSACF}
From \cite{Gamerman2006}, we give the estimation the autocorrelation function (ACF) at fixed lags and the approximated effective sample size (ESS) of the produced chain $\{x^{(n)}\}_{1}^{N_{s}}$. Define the autocovariance of lag $k\ (k\geq 0)$ of the chain as $\gamma_{k}=Cov_{\pi}(x^{(n)},x^{(n+k)})$, the autocorrelation function of lag $k$ as $\rho_{k}=\dfrac{\gamma_{k}}{\gamma_{0}}$ and  $$\tau=\sum_{k=-\infty}^{k=\infty} \rho_k=1+2\sum_{k=1}^{k=\infty} \rho_k.$$
If the series of autocorrelation is summable, $\tau$ is called  integrated autocorrelation time (IAT) as \cite{Green1992}. Consider that $x \in R^{d}$, we can derive the IAT $\tau_i, \ i=1,2,\cdots,d$ of each dimension $x_i$ similar as above.
Then, the ESS of each dimension can be given
$ESS_i=\dfrac{N_s}{\tau_i}$ and the monolithic ESS is defined by
\begin{equation}\label{eq:ESS}
ESS=\dfrac{N_s}{\tau_{max}},\ \ \tau_{max}=\max\limits_{i=1,2,\cdots,d} \{\tau_i\}
\end{equation}

\end{document}